\documentclass[11pt]{article}

\def\mathbb{\Bbb}
\usepackage{amsfonts,latexsym,amstext}
\usepackage{amsmath}
\usepackage{amssymb}

\usepackage{comment}

\newtheorem{theorem}{Theorem}[section]
\newtheorem{lemma}[theorem]{Lemma}
\newtheorem{proposition}[theorem]{Proposition}
\newtheorem{definition}{Definition}[section]

\newtheorem{remark}[theorem]{Remark}

\makeatletter
\@addtoreset{equation}{section}

\makeatother

\def\qed{{\hfill\hbox{\enspace${ \square}$}} \smallskip}
\def\sqr#1#2{{\vcenter{\vbox{\hrule height .#2pt \hbox{\vrule
 width .#2pt height#1pt \kern#1pt \vrule
width .#2pt} \hrule height .#2pt}}}}
\def\square{\mathchoice\sqr54\sqr54\sqr{4.1}3\sqr{3.5}3}

\def\ds{\begin{displaystyle}}
\def\eds{\end{displaystyle}}
\def\dis{\displaystyle }
\def\<{\langle }
\def\>{\rangle }

\def\R{\mathbb R}

\def\E{\mathbb E}
\def\P{\mathbb P}

\def\calb{{\cal B}}

\def\calf{{\cal F}}

\def\calh{{\cal H}}
\def\cali{{\cal I}}

\def\caln{{\cal N}}
\def\calo{{\cal O}}
\def\calp{{\cal P}}
\def\calt{{\cal T}}

\def\caly{{\cal Y}}

\title{Filtering of continuous-time Markov chains with noise-free observation
and applications}

\date{}
\author{Fulvia Confortola, Marco Fuhrman\\
Politecnico di Milano,
Dipartimento di Matematica\\
piazza Leonardo da Vinci 32, 20133 Milano, Italy\\
e-mail: fulvia.confortola@polimi.it, marco.fuhrman@polimi.it}

\begin{document}

\maketitle

\begin{abstract}
Let
$X$ be a continuous-time  Markov chain
in a finite set $I$, let $h$ be a mapping
of $I$ onto another set, and let $Y$
be defined by $Y_t=h(X_t)$, ($t\ge 0$).
We address the filtering problem for $X$ in terms of the
observation $Y$, which is not directly affected by noise.
We write down explicit equations for the filtering process
$
\Pi_t(i)=\P(X_t=i\,|\, \caly_t^0)$, ($i\in I, t\ge 0$),
where $(\caly_t^0)$ is the natural filtration of $Y$.
We show that $\Pi$ is a Markov process with the Feller property.
We also prove that it is a
piecewise-deterministic Markov process
in the sense of  Davis, and we identify its
characteristics explicitly. We finally solve an optimal
stopping problem for $X$ with partial observation,
i.e. where the moment of stopping is required to be a stopping time
 with respect to $(\caly_t^0)$.
\end{abstract}

\section{Introduction}

In the classical formulation of the filtering problem
in continuous time
the basic datum is a pair of stochastic processes $(X_t)_{t\ge 0}$
and $(Y_t)_{t\ge 0}$, defined on some probability space
$(\Omega,\calf,\P)$
with values in measurable spaces $(I,\cali )$ and
$(O,\calo )$ respectively. $X$ is called the unobserved (or signal) process
and $Y$ the observation process.
The filtering process is defined as
$$
\Pi_t(A)=\P(X_t\in A\,|\, \caly_t^0),\qquad A\in \cali, t\ge 0,
$$
where $\caly_t^0=\sigma(Y_s,\; s\in [0, t])$ are the $\sigma$-algebras
of  the filtration
generated by the observation process.
The filtering problem consists in a description of the measure-valued
process $\Pi$ and its properties. Often, $\Pi$ is shown to
satisfy some differential equations, called the filtering equations,
and in several cases it can be characterized as the unique solution
of such equations.

Various kinds of unobserved processes  have
been considered in the literature. $X$ is often taken to be a
diffusion process solution of a stochastic differential equations
driven by the Wiener process.
The case  of
 $X$ being a
Markov chain, or more generally a marked point process, is also
frequently addressed.

Concerning the observation process, the large majority
of cases considered in the literature are variants of the following
situation:
$Y$ takes values in the euclidean space $O=\R^m$ and has the form
\begin{equation}\label{osservaz}
  Y_t=\int_0^tH(X_s)\,ds+ \sigma W_t,\qquad t\ge 0,
\end{equation}
where $W$ is a standard Wiener process in  $\R^m$,
$H:I\to \R^m$ is a given function and
$\sigma$ is a constant,
$\sigma\neq 0$.

As is well known, special results are available in the linear gaussian
case, the so called Kalman-Bucy filter.

Recently, the following different model has been addressed by several
authors:
\begin{equation}\label{osservazperf}
  Y_t=h(X_t),\qquad t\ge 0,
\end{equation}
where $h:I\to O$ is a given function. We call $Y$
a noise-free observation, since
it is not directly affected by noise but rather
all  sources of randomness are included in the unobserved process
$X$ (in \cite{JoLGl}, $Y$  is
called perfect observation).

A basic
motivation for studying noise-free observations
arises in
connection with the following variant of
(\ref{osservaz})
\begin{equation}\label{osservazdipstato}
  Y_t=\int_0^tH(X_s)\,ds+ \int_0^t\sigma(X_s) \,dW_s,\qquad t\ge 0,
\end{equation}
where a
state-dependent diffusion coefficient $\sigma$ occurs. It can be proved
that this model can be converted into another one where the
observation process has two components: one is similar
to the traditional model (\ref{osservaz}) and the other
is noise-free, i.e., of the form (\ref{osservazperf}):
see \cite{TaAk}, \cite{KorRung}, \cite{CrKouXi}.

Noise-free observations have also been considered in connection
with the classical topic of the filter stability.
It has been discovered that several ergodicity properties of the filtering
process $\Pi$ fail to hold when the observation is noise-free:
see for instance
\cite{BaChiLip} and the discussion and references therein;
see also our remark \ref{instabfilter}.

Another motivation, pointed out in
\cite{JoLGl}, is that a clear picture of the noise-free
case may also lead to a better understanding of the limiting
behavior of the model (\ref{osservaz}) as $\sigma\to 0$.

Beside these motivations, it is our opinion that the case of
noise-free observation deserves attention in its own. The few existing results
are already mathematically interesting. More important, it seems to be
the natural mathematical model to describe
situations where the available observations
are indeed very accurate. For instance,
when randomness is introduced in order to represent
uncertainty about the state of an evolving dynamical system
it might be unnatural to introduce a noise affecting the observation
in order to fit the standard framework (\ref{osservaz})
unless this corresponds to a substantial description of
inaccuracy of measurements.

In spite of its interest, there are
few existing results on the case of noise-free observation,
with the important exception of the Kalman-Bucy filtering theory
(for the latter we limit ourselves to noting that
generalizations of (\ref{osservaz}) to the case of
degenerate noise acting on $Y$
date back at least at  the paper
\cite{BrJo}, and we
will not give detailed references for the linear gaussian case).
In fact, we are only aware of \cite{JoLGl} as the
 only paper entirely devoted to nonlinear filtering
 in the case of noise-free observation.
In \cite{JoLGl}, $X$ is defined as the solution
to a stochastic differential equation in $I=\R^n$, the observation takes
values in $O=\R^m$
and the function $h$ in (\ref{osservazperf}) is assumed
to satisfy special assumptions. The main result is that
the filtering equation can be formulated as a stochastic
equation on a submanifold of $\R^n$, and under appropriate conditions
the law of $\Pi_t$ admits a density with respect to the
surface measure.
These results are used in
\cite{CrKouXi} to study the model
(\ref{osservazdipstato}).

\medskip

The purpose of the present paper is a systematic study
of the filtering problem with noise-free
observation  when $X$ is assumed to be a time-homogeneous Markov chain
in a finite set $I$. Thus, our basic data will be
 a pair of finite sets $I$ and
$O$, a function $h:I\to O$ (assumed to be surjective
without loss of generality), the
rate
transition matrix $\Lambda$ of a Markov chain $X$ in $I$,
and  the noise-free
observation process defined by (\ref{osservazperf}). The filtering
process is specified by the finite set of scalar processes
$$
\Pi_t(i)=\P(X_t=i\,|\, \caly_t^0),\qquad i\in I, t\ge 0,
$$
where $\caly_t^0=\sigma(Y_s,\; s\in [0, t])$.

Our main results are the following. In section
\ref{filteringproblem}, after some notation and preliminary
results, we present the filtering equations: they are
a system of ordinary differential equations with jumps and
with random coefficients (depending on the observation process)
and  their unique solution provides a modification of the
filtering process $\Pi$. The method of proof is based on a
discrete approximation: we first write down the filtering equations
corresponding to observing the process $Y$ only at times $k2^{-n}$
($k=0,1,\ldots$) and we pass to the limit as $n\to \infty$.

The noise-free case has the following special feature:
if at some time $t$ we observe $Y_t=a\in O$ then we know that
$\Pi_t$ has support in the corresponding level set of the function
$h$, namely $h^{-1}(a)=\{i\in I\,:\, h(i)=a\}$. So the natural
state space of the process $\Pi$ consists of probability measures
on $I$ which are supported in one of the level sets $h^{-1}(a)$
($a\in O$). We call this space the effective simplex $\Delta_e$.
In section \ref{basicprop}, after introducing an appropriate
canonical set-up and solving the so called prediction
problem, we establish that
$\Pi$ is a Markov process in $\Delta_e$ with respect to
the filtration $(\caly_t^0)$:
see Proposition \ref{filtromarkoviano}.
We also recall a  known counterexample about the lack of ergodic properties
of the filtering process in the case of noise-free observation: see
Remark \ref{instabfilter}.

Since the trajectories of the observation process are piecewise constant,
the law of $Y$ is completely determined by the finite-dimensional distributions
of the process $Y_0, T_1, Y_{T_1},T_2, Y_{T_2},\ldots$
where $T_j$ denote the  jump times of $Y$.
In
section
\ref{sectionobservation} we find explicit formulae for
those distributions in terms of
the filtering process $\Pi$.
Although this is mainly a technical point in preparation
of the  results to follow, it  has some immediate
application: see for instance
{Proposition} \ref{exitime} where we
 prove an explicit formula for
the law of the exit time of a finite Markov chain from a given set,
a result that we could not find in the literature.

We note that in our model
new information
is available only at
jump times $T_j$.
Therefore it is not surprising that the filtering equations
prescribe a smooth, deterministic evolution of the trajectories
 $t\mapsto \Pi_t(\omega)$ among such jump times, and a jump of $\Pi$
 at each time $T_j$. An important class of Markov processes, having jumps
at some random times and otherwise evolving along a deterministic flow,
was introduced by M.H.A. Davis in \cite{Da2} and named
piecewise-deterministic Markov processes
 (PDPs).
In section
\ref{sectionpdp} we show that the
filtering process $\Pi$ is a PDP in the sense of Davis,
and we explicitly describe its characteristics
(the flow, the
 jump rate function and the transition measure):
 see
Theorem \ref{mainpdp}.
Since PDPs are  processes that have been extensively studied,
see for instance the book  \cite{Da}, a lot of known results
on PDPs
immediately applies to the filtering process. For instance,
a precise description of its generator is known, in terms of
the characteristics. Moreover, we are able to show the
Feller property for the process $\Pi$ using
arguments from \cite{Da}: see Proposition \ref{pdpfeller}, where
we prove in addition that the transition semigroup
of $\Pi$ is   strongly continuous in the space of continuous functions
on the state space
$\Delta_e$, equipped with the supremum norm.

In section
\ref{optimalstopping}
we study an optimal stopping problem for the Markov
chain $X$ with partial observation. The functional to be minimized
has classical form, but the moment of stopping   is subject to be
a stopping time with respect to the filtration $(\caly^0_t)$
generated by $Y$, i.e., it has to be based only on the observed process.
We follow the classical approach to first solve an optimal stopping
problem for the filtering process $\Pi$, appropriately formulated
and with complete observation, and then to show  how this gives
a solution to the original problem: see for instance
\cite{MSSZ} and the references therein
for this approach in a general framework.
Once more the theory of PDPs turns out to be a very useful tool here,
since the existence of an optimal stopping time
for  $\Pi$, as well as a characterization of the
value function and the stopping rule,
are a direct application of known results
on PDPs. We obtain corresponding results for the original
problem with partial observation.

\medskip

When $X$ is a finite Markov chain and the observation
has the classical form (\ref{osservaz}),
the corresponding process $\Pi$ is called the Wonham filter,
see \cite{Wo}.
In spite of the simple structure of the unobserved process $X$,
this case is still the subject of current study, see for instance
\cite{BaChiLip} for investigations on stability of the Wonham filter.
In
 section 3 of that paper the following
earlier example of noise-free observation
due to
\cite{DeZe}
is analyzed:
$X$ is  the Markov chain
in the   space $I=\{1,2,3,4\}$ with rate
transition matrix
\begin{equation}\label{lambdacontroes}
  \Lambda=\left(
\begin{array}{cccc}
-1&1&0&0
\\
0&-1&1&0
\\
0&0&
-1&1
\\
1&0&0&-1
\end{array}\right)
\end{equation}
and $Y$ takes values in $O=\{0,1\}$ with
the function $h$ in (\ref{osservazperf}) given by
\begin{equation}\label{hcontroes}
h(1)=h(3)=1, \quad h(2)=h(4)=0.
\end{equation}
In this specific case the filtering equations are deduced:
see Proposition 3.2 in \cite{BaChiLip}. The  method
of proof is different from ours and relies on
martingale methods.
It should be mentioned at this point that
several methods are known to prove that the filtering
process is a solution of the corresponding filtering equations,
 in the case of noisy observation
 (\ref{osservaz}). It is possible that some of them can be applied
to deduce the filtering equations in the general case of noise-free
observation as well. For instance one may
try to generalize Proposition 3.2 in \cite{BaChiLip} mentioned above,
or one may rely on the fact that the natural filtration of a jump process
(in our  case,  the filtration $(\caly_t^0)$) can be described
in a precise way: see for instance \cite{Da} or the earlier works
\cite{BoVaWo}. However in this paper we present an elementary proof
of the filtering equations,
based on discrete approximation, that is self-contained and does
not use deep results from the general theory of stochastic processes.

 We finally mention that the main results of this paper
 have been presented at the
  {\it First CIRM-HCM Joint Meeting: Stochastic Analysis, SPDEs,
Particle Systems, Optimal Transport}, Levico Terme (Italy),
January 24-30,
2010.

\section{The filtering problem}
\label{filteringproblem}

\subsection{Formulation}

The filtering problem will
be described starting from
the following basic objects, which are assumed to
be given throughout the paper.

\begin{enumerate}
\item $I$ is a finite set.

Elements of $I$, usually denoted by letters
$i,j,k\ldots$, are called states, and $I$ is called state space.

\item $\Lambda$ is a rate transition matrix on $I$ (sometimes called a $Q$-matrix,
see e.g. \cite{No}), i.e. a square matrix $(\lambda_{ij})_{i,j\in I}$ whose elements
are real numbers satisfying
$
\lambda_{ij}\ge 0$ for $ i\neq j$ and $\sum_{j\in I}\lambda_{ij}=0$ for all $i\in I$.

\item $h$ is a surjective function defined on $I$ taking values in another finite set $O$.

$O$ is called the observation space. Its elements will be denoted by letters
$a,b,c\ldots$.
The assumption that $h$ is surjective does not involve any loss of generality, since
$O$ may be replaced by the image of $h$ in all that follows.

\end{enumerate}

Suppose that on some
probability space $(\Omega,\calf,\P)$
a process $(X_t)_{t\ge0}$ is defined, taking values in $I$.
$(X_t)$ will be called the unobserved process.
We
assume that it is a Markov process with
generator $\Lambda$, i.e.
for every $t,s\ge0$ and every real function $f$ on $I$ we have
$$
\E [f(X_{t+s})|\calf^0_t]=(e^{s\Lambda}f)(X_t),\qquad \P-a.s.
$$
where
$
\calf^0_t=\sigma(X_s\,:\,s\in[0,t])
$, $t\ge0$,  are the $\sigma$-algebras of the natural filtration of $(X_t)$.
We denote  by $\mu$ the initial distribution of $(X_t)$, i.e.  $\mu(i)=\P(X_0=i)$, $i\in I$.

Next we define the observation process $(Y_t)_{t\ge0}$ by the formula
$$Y_t=h(X_t),\qquad t\ge 0,
$$
and  we introduce its
natural filtration  $(\caly^0_t)_{t\ge0}$, where
$$
\caly^0_t=\sigma(Y_s\,:\,s\in[0,t]),\qquad t\ge 0.
$$
The filtering problem consists in describing the conditional
distribution of $X_t$ given $\caly^0_t$, for all ${t\ge0}$. In
other words, we look for a description of the processes
$\P(X_t=i|\caly^0_t)$,  ${t\ge0}$, for all $i\in I$. The main
result of this section, Theorem \ref{teo-filtro}, states that an
appropriate modification of those processes, denoted by
$\Pi_t(i)$, are the unique solutions of suitable differential
equations, called the filtering equations, driven by the
observation process. In addition, these equations are explicitly
written and their solution is clearly described.

In order to present those equations we need to introduce some notation
and preliminary results, presented in the next subsections.

\subsection{The effective simplex and the flow}\label{simplexflow}

By $\Delta$ we denote the set of probability measures on $I$. $\Delta$ can be
naturally identified with the canonical simplex of $\R^N$, where $N$ is
the cardinality of $I$, i.e. with the set of  $\mu=(\mu(i))\in\R^N$
such that $\mu(i)\ge 0$ for all $i=1,\ldots,N$ and $ \sum_{i=1}^N\mu(i)=1$.

We note that the sets $h^{-1}(a)=\{i\in I\,:\, h(i)=a\}$ form a partition of $I$
as $a$ varies in $O$. We denote by $\Delta_a$ the set of probability measures on
$I$ supported in
$h^{-1}(a)$.  Each $\Delta_a$ can be considered as a subset
of $\Delta$, so an element in $\Delta_a$ is an $N$-dimensional vector
$\mu=(\mu(i))$ in $\Delta$ such that $\mu(i)=0$ if $h(i)\neq a$.
This way each $\Delta_a$ is a simplex in the euclidean space $\R^N$ and
it is a face
of $\Delta$, i.e. its vertices are also vertices of $\Delta$. Moreover,
as $a$ varies in $O$,
the vertices of the simplices $\Delta_a$ form a partition of the vertices of $\Delta$.

Finally we define $\Delta_e=\cup_{a\in O}\Delta_a$. This is a subset
of $\Delta$, which we call the effective simplex. It is a proper subset unless
$h$ is constant. $\Delta_e$ is obviously a compact space.

As a general rule probability measures $\mu$ on $I$, or equivalently elements of $\Delta$,
will be considered as row vectors $(\mu(i))$ of dimension $N$ equal to the cardinality of $I$,
whereas functions $f:I\to\R$ will be identified with the $N$-dimensional column vector
$(f(i))$ consisting of its values. The integral of $f$ with respect to $\mu$ is denoted
simply $\mu f=\sum_{i\in I}f(i)\mu(i)$.
Given a row vector $(\nu(i))$
(not necessarily a probability measure)
and a column vector  $(f(i))$,
 both real and $N$-dimensional,  we denote
by $f\ast \nu$ the row vector obtained by pointwise multiplication, i.e.
having components
$(f\ast \nu)(i) =f(i)\nu(i)$, $i\in I$.

In describing properties of the filtering process a basic role will be played
by a flow $\phi$ on the effective simplex $\Delta_e$. The following lemma is used
to define $\phi$ by means of a differential equation.

\begin{proposition}\label{lemmaflow} For every $a\in O$ and $x\in\Delta_a$
the differential equation
\begin{equation}\label{eqflow}
y'(t)= 1_{h^{-1}(a)}* (y(t)\Lambda) - (y(t)\Lambda
1_{h^{-1}(a)})\,y(t),\qquad t\ge 0,
\end{equation}
with initial condition $y(0)=x$,
has a unique global  solution
 $y:\R_+\to\R^N$, where $N$ denotes  the cardinality of $I$.

 Moreover $y(t)\in\Delta_a$ for all $t\ge 0$.
\end{proposition}

We will write $\phi_a(t,x)$ instead of $y(t)$, to stress dependence on $a$ and $x$.
By standard results on ordinary differential equations,
$\phi_a$ is a continuous function of
$(x,t)$ and it enjoys the flow property
$\phi_a(t,\phi_a(s,x))=\phi_a(t+s,x)$, for $t,s\ge0$, $a\in O$ and $x\in \Delta_a$.
$\phi_a$ is the flow associated to the vector field
$$
F_a(y)=1_{h^{-1}(a)}* (y\Lambda)
- (y\Lambda 1_{h^{-1}(a)})\,y,\qquad y\in \Delta_a,
$$
on $\Delta_a$.
To simplify the notation a little, it is convenient to define a global flow $\phi$ on $\Delta_e
=\cup_{a\in O}\Delta_a$ in the obvious way, setting
$\phi(t,x)=\phi_a(t,x)$ if $x\in \Delta_a$. This way $\phi(t,\cdot)$ is a function
$\Delta_e\to\Delta_e$
leaving each set $\Delta_a$ invariant.

\medskip

The rest of this subsection is devoted to
the proof of
Proposition \ref{lemmaflow}
by means of a suitable  viability theorem. For this we  have to recall
the notion of contingent cone (see, e.g.,  \cite{Au} Chapter 1).
\begin{definition}
Let $X$ be a normed space, $K$ be a nonempty subset of $X$
and $x$ belong to $K$. The contingent cone to $K$ at $x$ is the set
$$T_K(x)=\left\{ v \in X \,:\, \liminf_{h \rightarrow 0^+} \frac{d_K(x+hv)}{h}=0 \right\},
$$
where $d_K(y)$ denotes the distance of $y$ to $K$, defined by
$$d_K(y):= \inf_{z \in K} ||y-z||.$$

In other words, $v$ belongs to $T_K(x)$ if and only if there exist
a sequence of $h_n >0$ converging to $0+$ and a sequence of $v_n
\in X$ converging to $v$ such that $$\forall n \geq 0, \,\, x +h_n
v_n \in K$$

\end{definition}

We need to compute the contingent cone to the set $ \Delta_a$.

\begin{lemma} The contingent cone
$ T_{\Delta_a}(y)$ to $\Delta_a$ at $y \in\Delta_a$ is the cone of elements $v \in \mathbb{R}^N$ satisfying
\begin{equation}\label{contingent cone}
\sum_{i \in h^{-1}(a)} v_i=0 \quad  \&\quad \left\{
 \begin{array}{ll}
  v_i = 0, &   \hbox{ for } i \notin h^{-1}(a)\\
  v_i \geq 0, & \hbox{ if }  y_i=0  \hbox{ for } i \in h^{-1}(a)                                                                                                          \end{array}
   \right.
\end{equation}
\end{lemma}
\noindent {\bf Proof.} Let us take $v \in T_{\Delta_a}(y)$. There
exist sequences $h_n>0^+$ converging to $0$ and $v_n$ converging
to $v$ such that $z_n:=y+h_n v_n$ belongs to $\Delta_a$ for any $n
\geq 0$. Then
$$\sum_{i \in h^{-1}(a)} v_{n_i}= \frac{1}{h_{n}}\left(\sum_{i \in h^{-1}(a)} z_{n_i} -\sum_{i \in h^{-1}(a)} y_i\right)=0,$$
so that $\sum_{i \in h^{-1}(a)} v_{n_i}=0$. On the other hand, if
$i \notin h^{-1}(a)$, then $y_i=z_{n_i}=0$, then $v_{n_i}=
\frac{z_{n_i}-y_i}{h_{n}}=0$, so that $v_{n_i} =0$.  If $y_i=0$
for $i \in h^{-1}(a)$, then $v_{n_i}= \frac{z_{n_i}}{h_{n}}\geq0$,
since $z_n$ belongs to $\Delta_a$. So $v_{n_i} \geq 0$.

Conversely, let us take $v$ satisfying (\ref{contingent cone})
and deduce that $z:= y+tv$ belongs to $\Delta_{a}$ for $t$ small enough. First we have
$$\sum_{i \in h^{-1}(a)} z_i=\sum_{i \in h^{-1}(a)} y_i +t \sum_{i \in h^{-1}(a)} v_i=\sum_{i \in h^{-1}(a)} y_i=1.$$
Second, if $i \notin h^{-1}(a)$, then $v_i=0$ and since $y \in \Delta_a$ also $y_i=0$. Hence if $i \notin h^{-1}(a)$ then $z_i=0$. Moreover if $i \in h^{-1}(a)$, then $z_i \geq0$. In fact if $y_i=0$ then $
z_i=tv_i \geq 0$, because in this case $v_i$ is nonnegative.
If $y_i >0$ it is sufficient to take $t \leq \frac{y_i}{|v_i|}$ for having $z_i \geq 0$.
Hence $z$ belongs to $\Delta_a$.
\qed

\medskip

\noindent {\bf Proof of proposition \ref{lemmaflow}.} $\Delta_a$
is closed and is a viability domain of $F_a$. This means that, for
all $y \in \Delta_a$, $F_a(y)$ belongs to $T_{\Delta_a}(y)$, the
contingent cone to $\Delta_a$ at $y$. In fact $v=F_a(y)$ satisfies
(\ref{contingent cone}):
\begin{eqnarray*}\sum_{i \in h^{-1}(a)} v_i & =& \sum_{i\in h^{-1}(a)}
[1_{h^{-1}(a)}* (y\Lambda)]_i - \sum_{i\in h^{-1}(a)}(y\Lambda 1_{h^{-1}(a)})y_i\\
 & = &[1_{h^{-1}(a)}* (y\Lambda)]1_{h^{-1}(a)} -
(y\Lambda 1_{h^{-1}(a)})y1_{h^{-1}(a)}\\
& = & (y\Lambda 1_{h^{-1}(a)})-
(y\Lambda 1_{h^{-1}(a)})y1_{h^{-1}(a)}\\
 & = & (y \Lambda
1_{h^{-1}(a)})(1-y 1_{h^{-1}(a)})=0,
\end{eqnarray*}
since $y \in \Delta_a$ and $\sum_{i \in h^{-1}(a)}y_i=y
1_{h^{-1}(a)}=1$. Moreover if $ i \notin h^{-1}(a)$, then
$v_i=[1_{h^{-1}(a)}* (y\Lambda)]_i - (y\Lambda
1_{h^{-1}(a)})y_i=0$, since the components of $v_i$ are both equal
to zero. If $y_i=0$ for $i \in h^{-1}(a)$, then
$$v_i=[1_{h^{-1}(a)}* (y\Lambda)]_i= \sum_{j\in h^{-1}(a),j
\neq i } y_j \lambda_{ji} \geq 0,$$ since $\lambda_{ji} \geq 0$
for $j \neq i$ and $y_j \geq0$ for $j\in h^{-1}(a)$.

Moreover  $|F_a(y)| \leq C(1+|y|)$ for all $y \in
\Delta_a$. Then by Theorem 1.2.4 in \cite{Au},
page 28, $\Delta_a$ is viable under $F_a$: for every initial state
$x$ there exists a solution $y(\cdot)$ to differential equation
(\ref{eqflow}) such that $y(t) \in \Delta_a$ for each $t \in [0,
\infty)$. The uniqueness follows from the fact that the function
$F_a(y)$ is locally Lipschitz on $\mathbb{R}^N$.  \qed

\subsection{The operator $H$}
\label{sectionoperatorh}

For every $a\in O$ we define a function $H_a$, mapping
row vectors $\mu\in\R^N$ to row vectors $H_a[\mu]\in\R^N$ defined
for $i\in I$ by
$$
H_a[\mu](i)=\left\{
\begin{array}{ll}
0,
&{\rm if\;} h(i)\neq a,
\\\dis
\frac{\mu(i)}
{\sum_{j\in h^{-1}(a)}\mu(j)},
&{\rm if\;} h(i)= a, \;\sum_{j\in h^{-1}(a)}\mu(j)\neq 0,
\\
\nu_a
&{\rm if\;} h(i)= a, \;\sum_{j\in h^{-1}(a)}\mu(j)= 0,
\end{array}
\right.
$$
where $\nu_a$ is a fixed arbitrary probability on $I$ supported in $h^{-1}(a)$,
whose exact values are irrelevant.
Using the notation introduced in the previous section we may write
$$
H_a[\mu]=\left\{
\begin{array}{ll}\dis
\frac{1}{\mu 1_{h^{-1}(a)}}1_{h^{-1}(a)}\ast \mu,
&{\rm if\;} \mu 1_{h^{-1}(a)}\neq0,
\\
\nu_a
&{\rm if\;} \mu 1_{h^{-1}(a)}=0.
\end{array}
\right.
$$
We note that if $\mu(i)\ge 0$ for all $i$ then $H_a[\mu]$ is a probability
measure on $I$ supported on $h^{-1}(a)$, i.e. an element of $\Delta_a$.
If in addition $\mu$ is a probability then $H_a[\mu]$ is the corresponding conditional
probability given the event  $\{h=a\}$.

The relevance of the operator $H$ to the filtering problem is well
known. Indeed suppose that  $(\bar{X}_{k})_{k\in \mathbb{N}}$ is a
discrete-time Markov process in $I$ with transition matrix $P$ and
initial distribution $\mu$ and the observation process is defined
by $\bar{Y}_{k}=h(\bar{X}_{k})$. Then, see for instance
\cite{CMR}, the discrete-time filtering process defined for  $i
\in I$ by the formula
$\bar{\Pi}_k(i):=\P(\bar{X}_{k}=i|\bar{Y}_0,\ldots,\bar{Y}_{k})$
satisfies the recursive equations
\begin{equation}\label{filtrodiscreto}
\left\{\begin{array}{l}
\bar{\Pi}_{0}= H_{\bar{Y}_{0}}[\mu],
\\
\bar{\Pi}_{k}= H_{\bar{Y}_{k}}[
\bar{\Pi}_{k-1}P], \qquad k\ge 1.
\end{array}
\right.
\end{equation}

\subsection{The filtering equation}
Let $(T_j)_{j\ge 1}$ denote the  sequence of  jumps times of
$(Y_t)$, with the convention that $T_j=\infty$ for all $j$ if no
jump occurs and $T_{j}<T_{j+1}=T_{j+2}=\ldots=\infty$ if precisely
$j$ jumps occur. We set $T_0=0$.

We define a process
  $(\Pi_t)$  and we will eventually prove that it is a modification
  of the filtering process.
For every $\omega \in \Omega$ we consider the corresponding
trajectory $Y_t(\omega)$ and jumps times $T_j(\omega)$. Next we set
$\Pi_0(\omega)=  H_{Y_0(\omega)}[\mu]$ and
 for $j\ge1$,
 \begin{equation}\label{eq-pifiltro-equiv}\begin{array}{ll}
\Pi_t(\omega)=\phi(t-T_{j-1}(\omega), \Pi_{T_{j-1}(\omega)})& {\rm    \;for\;}  T_{j-1}(\omega)\le
t <T_j(\omega),
\\
\Pi_{T_j-}(\omega)=\phi(T_j(\omega) -T_{j-1}(\omega), \Pi_{T_{j-1}(\omega)}),& {\rm    \;if\;}
T_j(\omega)< \infty,
\\
\Pi_{T_j}(\omega)=H_{Y_{T_j}(\omega)}[\Pi_{T_j-}(\omega)\Lambda],& {\rm    \;if\;}  T_j(\omega)<
\infty,
\end{array}
\end{equation}
where $\phi $ is the global flow defined in  Section
\ref{simplexflow}.

Note that $\Pi_{T_j-}$ and $\Pi_{T_j}$ are only defined on $\{T_j<\infty\}$ and on this
set $\Pi_{T_j-}(\omega)$ is the usual left limit $\lim_{t\to {T_j}(\omega),
t<{T_j}(\omega)}\Pi_t(\omega)$.

It follows from the properties of the operator $H$ that $\Pi_{T_j}$ is supported
on $h^{-1}(Y_{T_j})$, $(j\ge 0)$; equivalently,
 $\Pi_{T_j}\in\Delta_{Y(T_j)}$. Moreover, since the flow $\phi$ leaves each set
$\Delta_a$ invariant, we conclude that
 $\Pi_{t}$ is supported
on $h^{-1}(Y_{T_{j-1}})$ for $T_{j-1}(\omega)\le
t <T_j(\omega)$,  $(j\ge 1)$.

The process
  $(\Pi_t)$ can also be described in the following way. Its starting point is
$\Pi_0=  H_{Y_0}[\mu]$ and it has jumps precisely at times $T_j$
 $(j\ge1)$. At each jump time it jumps from $\Pi_{T_j-}$ to
 $\Pi_{T_j}=H_{Y_{T_j}}[\Pi_{T_j-}\Lambda]$.  Among jump times,
 trajectories evolve deterministically. Namely, for
 $ T_{j-1}\le
t <T_j$, we have $\Pi_t(i)=0$ if $i\notin h^{-1}(Y_{T_{j-1}})$ and
\begin{equation}\label{eqdiffintervallo}
     \Pi_t^{'}(i)=
 (\Pi_t\Lambda)(i) -
(\Pi_t\Lambda 1_{h^{-1}(Y_{T_{j-1}})})\Pi_t(i)
\end{equation}
if $i\in h^{-1}(Y_{T_{j-1}})$,  where the time derivative is understood as the right derivative
at time $t=T_{j-1}$. Using previously introduced notation we can write
the differential equation
in vector form
$$
\Pi_t^{'}=
1_{h^{-1}(Y_{T_{j-1}})}* (\Pi_t\Lambda) -
(\Pi_t\Lambda 1_{h^{-1}(Y_{T_{j-1}})})\Pi_t \quad
\quad  t \in [T_{j-1},T_j)
$$
and we could even replace $Y_{T_{j-1}}$ by $Y_{t}$, since
$Y$ is constant among jump times.

Finally, it is also possible to describe the process $(\Pi_t)$ by a single
integral equation: for every  $\omega\in\Omega$,
\begin{equation}\label{eq-pifiltro}
\begin{array}{lll}
\Pi_t(\omega)&=&\dis  H_{Y_0(\omega)}[\mu]+\int_0^t\left\{
1_{h^{-1}(Y_s(\omega))}* (\Pi_s(\omega)\Lambda) -
(\Pi_s(\omega)\Lambda
1_{h^{-1}(Y_s(\omega))})\Pi_s(\omega)\right\}ds
\\&&\dis
+\sum_{0<T_j(\omega)\le t}\left\{
H_{Y_{T_j}(\omega)}[\Pi_{T_j-}(\omega)\Lambda]-\Pi_{T_j-}(\omega)\right\},
\qquad t\ge0.
\end{array}
\end{equation}
Note that only  indices $j\ge 1$ enter the sum, since
$0=T_0(\omega)<T_j(\omega)$.

We will refer to  (\ref{eq-pifiltro}) as the filtering equation. This is justified by next
theorem.

\begin{theorem}\label{teo-filtro}
 The process
 $(\Pi_t)$, defined by equation (\ref{eq-pifiltro}), is a modification
 of the filtering process: for  ${t\ge0}$ and $i \in I$ we have
$\Pi_t(i)=\P(X_t =i|\caly^0_t)$ $\P$-a.s.
 \end{theorem}

Recall that in the definition of the operator $H$
we used some arbitrary probabilities, denoted $\nu_a$. It is worth
noting that they do not affect the definition of the filtering process, in the
sense that we have, $\P$-a.s.,
\begin{equation}\label{dendiversodazero}\Pi_{T_j}=
H_{Y_{T_j}}[\Pi_{T_j-}\Lambda]=
\frac{\Pi_{T_j-}\Lambda}{\Pi_{T_j-}\Lambda
1_{h^{-1}(Y_{T_j})}}
\end{equation}
for every $j\ge 1$. Therefore, the filtering equation
(\ref{eq-pifiltro}) can also be written as follows:
$$
\begin{array}{lll}
\Pi_t(\omega)&=&\dis  H_{Y_0(\omega)}[\mu]+\int_0^t\left\{
1_{h^{-1}(Y_s(\omega))}* (\Pi_s(\omega)\Lambda) -
(\Pi_s(\omega)\Lambda
1_{h^{-1}(Y_s(\omega))})\Pi_s(\omega)\right\}ds
\\&&\dis
+\sum_{0<T_j(\omega)\le t}\left\{
\frac{\Pi_{T_j-}(\omega)\Lambda}{\Pi_{T_j-}(\omega)\Lambda
1_{h^{-1}(Y_{T_j}(\omega))}}
-\Pi_{T_j-}(\omega)\right\},
\qquad t\ge0.
\end{array}
$$
To prove (\ref{dendiversodazero}) it is enough to show that the denominator
$\Pi_{T_j-}\Lambda
1_{Y_{T_j}}$ never vanishes. This is the content of
 lemma \ref{P-Tj}, which is stated in a slightly more general form,
as required in the proof of
theorem \ref{teo-filtro}
that will follow.

\begin{lemma} \label{P-Tj}Let $(\Pi_t)_{t \geq 0}$ be a  process with nonnegative
components $\Pi_t(\omega,i)$ satisfying,  for some fixed
 $j \geq 1 $, the following conditions:
 \begin{enumerate}
 \item
 for  $ T_{j-1}(\omega)<\infty$ and
 $ T_{j-1}(\omega)\le
t <T_j(\omega)$, we have $\Pi_t(\omega,i)=0$ if $i\notin h^{-1}(Y_{T_{j-1}}(\omega))$ and
\begin{equation}\label{P-eq}
 \Pi_t^{'}(\omega,i)=
 (\Pi_t(\omega)\Lambda)(i) -
(\Pi_t(\omega)\Lambda 1_{h^{-1}(Y_{T_{j-1}}(\omega))})\Pi_t(\omega ,i)
\end{equation}
if $i\in h^{-1}(Y_{T_{j-1}}(\omega))$;

\item for every
$t \geq 0$
\begin{equation}\label{Pnfiltro}\Pi_t (i) 1_{T_{j-1}\leq t <T_j}
=\P(X_t=i|\mathcal{Y}^0_t)1_{T_{j-1}\leq t <T_j} \qquad \P-{a.s. }
\end{equation}
\end{enumerate}
Then \begin{equation}\label{Pi-Tj}\Pi_{T_j-} \Lambda
1_{h^{-1}(Y_{T_j})}>0 \qquad \P-{ a.s.\; on\;} \{T_j<\infty\}.
\end{equation}
\end{lemma}
\noindent {\bf Proof.} First we prove that for all $ i\in
I$
\begin{equation}\label{*}
\{  \Pi_{T_j-}(i)=0,T_j<\infty\} \subset \{ X_{T_j-} \neq i,T_j<\infty
\}, \qquad \P-a.s.
\end{equation}

We fix $\omega \in \Omega$ such that
$T_{j-1}(\omega)<\infty$,
and argue for $T_{j-1}(\omega) \leq t <
T_j(\omega)$. Define
$$E_t:= \exp \left( \int_0^t \Pi_s \Lambda
1_{h^{-1}(Y_{T_{j-1}})}ds\right),\qquad
U_t:=E_t \Pi_t.
$$
Since
$E_t>0$, $U_t$  and  $\Pi_t$ have the same support.
Then $U_t(i)=0$ if $i\notin h^{-1}(Y_{T_{j-1}})$. Now let
$i\in h^{-1}(Y_{T_{j-1}})$.
It follows
from  (\ref{P-eq}) that $U$
satisfies
$U^{'}_t(i) = (U_t \Lambda)(i)$
so that
$$U^{'}_t(i)= \sum_{k\in h^{-1}(Y_
{T_{j-1}})} U_t(k) \lambda_{ki}\geq  U_t(i)\lambda_{ii}.$$ This
implies that $U_t(i) \geq e^{ \lambda_{ii}(t-s)}U_{s}(i)$
for $T_{j-1}\le s\le t<T_j$, hence  $U_{s}(i)>0$ implies $U_t(i)>0$ for all $t \in
[s,T_j)$ and in particular
 $U_{T_{j-1}}(i)>0$ implies $U_t(i)>0$ for all $t \in
[T_{j-1},T_j)$.
The same result holds for
$\Pi_t(i)$, so we conclude
that, for every $i\in I$,
the process
\begin{equation}\label{processoduesalti}
    1_{\Pi_t(i)=0} 1_{T_{j-1} \leq t < T_j}
\end{equation}
has at most two jumps on $[0,\infty)$.
We also deduce that, if $T_j<\infty$,
 $$\Pi_{T_j-}(i)= \lim_{t\rightarrow T_j,
t< T_j} \Pi_{t }(i)= \lim_{t\rightarrow T_j, t< T_j} E_{t}^{-1}
U_{t} \geq  E_{T_j}^{-1} e^{\lambda_{ii}(T_{j} -s)}U_{s}(i),
\qquad
T_{j-1}\le s <T_j$$
and we also conclude that,
for every $i\in I$, on $\{T_j<\infty\}$,
\begin{equation}\label{pisiannulla}
\Pi_{T_j-}(i)=0 \qquad\Longleftrightarrow \qquad \Pi_{t}(i)=0 \; \forall t \in [T_{j-1},T_j).
\end{equation}
On the other hand, by  (\ref{Pnfiltro}), for every $t\ge 0$,
\begin{eqnarray*}&&
\P(X_t=i,\Pi_t(i)=0, T_{j-1}  \leq  t < T_j)= \E[1_{X_t=i} 1_{\Pi_t(i)=0} 1_{ T_{j-1} \leq t < T_j}]\\
&&=\E[ 1_{\Pi_t(i)=0} 1_{ T_{j-1} \leq t <
T_j}\P(X_t=i|\mathcal{Y}_t)]= \E[1_{\Pi_t(i)=0} 1_{ T_{j-1} \leq t
< T_j}\Pi_t(i)]=0.
\end{eqnarray*} Therefore
 the process
$1_{X_t=i} 1_{\Pi_t(i)=0} 1_{T_{j-1} \leq t < T_j}$
is a modification of the zero process. Recalling the process
in (\ref{processoduesalti}) and the fact that $X$ is piecewise constant, it is also
 indistinguishable from zero.

\noindent
Now we are ready to prove (\ref{*}), which is equivalent to
$\P(\Pi_{T_j-}(i)=0,
X_{T_j-}=i,T_j<\infty )=0$.
Suppose that for some $\omega \in \Omega$ we have $T_j<\infty$,
$\Pi_{T_j-}(i)=0$ and $
X_{T_j-}=i$. By (\ref{pisiannulla}) we also have
$\Pi_{t}(i)=0$ for all $t \in [T_{j-1},T_j)$ and denoting by
$S\ge T_{j-1}$ is the last jump time of
the chain $X$ before $T_j$ it follows that
$\Pi_{t}(i)=0$ and $X_t=i$ for all $t \in [S,T_j)$.
However, since $\P$-a.s.
we have $1_{X_t=i} 1_{\Pi_t(i)=0} 1_{T_{j-1} \leq t < T_j}=0$
for every $t$, this is only possible
with  zero probability.

Now we are able to prove  (\ref{Pi-Tj}). Noting that
$\Pi_{T_{j}-}$ is supported on
$h^{-1}(Y_{T_{j-1}})$, we
have, on $\{T_j<\infty\}$,
$$\Pi_{T_j-} \Lambda 1_{h^{-1}(Y_{T_j})}=
\sum_{i \in h^{-1}(Y_{T_{j-1}}),k \in h^{-1}(Y_{T_{j}})} \Pi_{T_j-}(i) \lambda(i,k).
$$
Since $T_j$ is a jump time for $Y$, the sets
$h^{-1}(Y_{T_{j-1}})$ and $h^{-1}(Y_{T_{j}})$ are disjoint, so that
 in the previous sum all nonzero terms correspond to indices $i\neq k$
 and consequently all terms are nonnegative.
 Since $X_{T_j}\in
h^{-1}(Y_{T_j})$ it follows that
$$\Pi_{T_j-} \Lambda 1_{h^{-1}(Y_{T_j})}\ge
\sum_{i \in h^{-1}(Y_{T_{j-1}})} \Pi_{T_j-}(i) \lambda(i,X_{T_j})
= \sum_{i \in I} \Pi_{T_j-}(i) \lambda(i,X_{T_j})
\qquad {\rm on\;} \{T_j<\infty\}.
$$
To conclude the proof it is enough to show that
 $\P(\sum_{i \in I} \Pi_{T_j-}(i) \lambda(i,X_{T_j})=0,T_j<\infty)=0$.
 Since again all terms in the last sum are nonnegative we have
$$\begin{array}{lll}
\P\left(\sum_{i \in I} \Pi_{T_j-}(i) \lambda(i,X_{T_j})=0\right)
&= &\P\left( \cap_{ i \in I} [\Pi_{T_j-}(i) \lambda(i,X_{T_j})=0,T_j<\infty]\right)  \\
& =& \P\left( \cap_{ i \in I} [\{\Pi_{T_j-}(i)=0,T_j<\infty\}
\cup \{ \lambda(i,X_{T_j})=0,T_j<\infty\}]\right)\\
& \leq & \P\left( \cap_{ i \in I} [\{X_{T_j-} \neq i,T_j<\infty\}
\cup \{ \lambda(i,X_{T_j})=0,T_j<\infty\}]\right)  \\
& =& \P( \cap_{ i \in I}[ \{X_{T_j-} \neq i,  \lambda(i,X_{T_j})>0,T_j<\infty \} \\
 & & \dis{\qquad  \qquad  \qquad \, \cup \{ \lambda(i,X_{T_j})=0,T_j<\infty\}])}
\end{array}
$$
where we have used  (\ref{*}) in the inequality. But clearly
$$
\P( \cap_{ i \in I}[ \{X_{T_j-} \neq i,  \lambda(i,X_{T_j})>0,T_j<\infty \}
  \cup \{ \lambda(i,X_{T_j})=0,T_j<\infty\}])=0
$$
since,  at  time
$T_j$, the chain $X$ must jump from some state $i \in I$ to  $X_{T_j}$. \qed

\medskip

The rest of this section is devoted to the proof of
 Theorem \ref{teo-filtro}.
The method we use is to construct a sequence $(\Pi^n_t(i))_{t\ge 0,i\in I}$
of approximate filtering processes, each corresponding to
observing  the process $Y$ in
a discrete set $\calp_n$  of times. $\Pi^n_t(i)$
are constructed so as to converge to $\P(X_t=i|\caly^0_t)$
as $n\to\infty$.
Then we write down explicit filtering equations for $\Pi^n_t(i)$.
Passing to the limit in
these equations we
prove that $\Pi^n_t(i)$ converges to the solution
$\Pi_t(i)$ of  (\ref{eq-pifiltro}). This way we
identify
$\Pi_t(i)$ with a modification of the filtering process
$\P(X_t=i|\caly^0_t)$.

\medskip

For all  integers $n,k\ge 0$ we
set $t^n_k= 2^{-n} k$ and consider the grid
 $\calp_n:=\{t^n_k\}_{k \ge 0}$ with mesh
$t^n_{k+1}-t^n_k=2^{-n}$.
For fixed $n$ and for $ t \geq 0$
we introduce the
$\sigma$-algebras $\caly^{n}_t=\sigma(Y_{t_0^n},\ldots,
Y_{t_k^n}: \, t_k^n \leq t)$.
The filtration $(\caly^{n}_t)_{t\ge 0}$ corresponds to observing the process
$Y$ only at times $t_k^n$.

For any $t\ge 0$, we have
$\caly^{n}_t \subseteq \caly^{n+1}_t$ and
$\caly^0_t$ is generated by $\cup_{n } \caly^{n}_t$, so
by a martingale convergence theorem we have
\begin{equation}\label{convfiltrovero}
    \lim_{n \rightarrow
\infty}\P(X_t=i|\caly^n_t)=\P(X_t=i|\caly^0_t),\qquad
\P \mbox{ -a.s.},\;i\in I.
\end{equation}
Next
we show that
the filtering
processes $\P(X_t=i|\caly^n_t)$ have  modifications, denoted $\Pi_t^n(i)$,
which satisfy explicit filtering equations.
Recalling that $\mu$ denotes the initial distribution of $X$, let us
 introduce processes $(\Pi^n_t(i))_{t\ge 0,i\in I}$ as follows:
for all $n\ge 0$ and $\omega\in\Omega$ define
\begin{eqnarray}
                               \Pi_0^n(\omega) &=& H_{Y_0(\omega)}[\mu], \label{Pn1}\\
    \Pi_t^n(\omega) &=&  \Pi_{t^n_{k-1}}^n(\omega)e^{(t-t_{k-1}^n) \Lambda }
     \; \hbox{ if }\; t_{k-1}^n\le t <t_k^n , \label{Pn2}\\
     \Pi_{t_k^n}^n(\omega) &=&  H_{Y_{t_k^n}(\omega)}[\Pi_{t_k^n-}^n(\omega)]
     \label{Pn3}
     \end{eqnarray}
for $k\ge 1$, where of course
   $ \Pi_{t_{k}^n-}^n(\omega) =  \Pi_{t_{k-1}^n}^n(\omega)e^{(t_{k}^n-t_{k-1}^n) \Lambda }
   = \Pi_{t_{k-1}^n}^n(\omega)e^{2^{-n}\Lambda }$.

\begin{lemma} \label{proprietaPn}
For  ${t\ge0}$ and $i \in I$ we have
$\Pi_t^n(i)=\P(X_t =i|\caly^n_t)$ $\P$-a.s.
\end{lemma}

 \noindent {\bf Proof.}  We fix $n$ and to
simplify the notation  we write $t_k$
instead of $t_k^n$.
 For every $n$ and $t$ we denote $\bar{\Pi}^n_t(i)$ an arbitrary version of
 $\P(X_t=i| \caly^n_t)$.
  We consider the discrete-time processes $(X_{t_k})_{k\ge 0}$
 and $(Y_{t_k})_{k\ge 0}$ obtained evaluating
$(X_t)_{t \geq 0}$ and $(Y_t)_{t \geq 0}$
at  times $t_k$.  $(X_{t_k})_{k\ge 0}$ is a Markov chain with initial law $\mu$ and
transition matrix $P= e^{(t_k-t_{k-1})\Lambda}=e^{2^{-n}\Lambda}$.
Therefore, for every $k\ge 0$,
$\bar{\Pi}_{t_k}^n(i)=\P(X_{t_k}=i|Y_0,\ldots,Y_{t_k})$
satisfies, $\P$-a.s.,
 \begin{equation}\label{ricors1}
\bar{\Pi}_0^n=H_{Y_0}[\mu]
 \end{equation}
 and the recursive equations
 \begin{equation}\label{ricors2}
    \bar{\Pi}_{t_k}^n= H_{Y_{t_k}}[
\bar{\Pi}_{t_{k-1}}^nP]=
H_{Y_{t_{k}}}[\bar{\Pi}_{t_{k-1}}^ne^{(t_k-t_{k-1})\Lambda}],
 \end{equation}
compare (\ref{filtrodiscreto}).

Next
 we fix $ t_{k-1}\le t <t_{k}$, we take $f:I\to \R$ (identified with
 $f\in\R^N$) and, using the fact that
 $ \sigma(Y_{t_j}) \subset
\sigma(X_{t_j})$ for every $j$, and the Markov property of $(X_{t})_{t\ge 0}$,
we obtain, $\P$-a.s.,
\begin{eqnarray*}\E(f(X_t)| \caly^n_t) & = &\E(f(X_t)|Y_0, \ldots, Y_{t_{k-1}})
\\
&=&\E[ \E(f(X_t)|X_0, \ldots, X_{t_{k-1}})| Y_0, \ldots, Y_{t_{k-1}}]\\
 & = &\E[(e^{(t-t_{k-1}) \Lambda} f)(X_{t_{k-1}})|Y_0, \ldots, Y_{t_{k-1}}]
 \\
 &=&\bar{\Pi}_{t_{k-1}}^n e^{(t-t_{k-1}) \Lambda } f,
 \end{eqnarray*}
which shows that
  \begin{equation}\label{ricors3}
\bar{\Pi}_t^n= \bar{\Pi}_{t_{k-1}}^ne^{(t-t_{k-1}) \Lambda }.
 \end{equation}

Comparing  (\ref{ricors1}),
(\ref{ricors2}),  (\ref{ricors3}) with
(\ref{Pn1}),
 (\ref{Pn2}), (\ref{Pn3}), we conclude that  $\Pi^n$ and $\bar{\Pi}^n$ are
 modifications of one another, and this proves the lemma.
\qed

Equations (\ref{Pn1})-(\ref{Pn3}) describe the time evolution of the
trajectories of $\Pi_t^n$.  For $k \geq 1$,
$\Pi_t^n$ satisfies the differential equation
\begin{equation}\label{eq-fil-app}(\Pi_t^n)'= \Pi_t^n \Lambda
\quad \mbox{     for } t_{k-1}^n\le t <t_{k}^n.
\end{equation}
where the derivative is understood as the
right derivative at $t=t^n_{k-1}$. At each time
$t_k^n$, $\Pi_t^n$ takes on the value
$H_{Y_{t^n_{k}}}[\Pi_{t^n_{k}-}^n]$, thus possibly making a jump.
Therefore equations (\ref{Pn1})-(\ref{Pn3}) can also be written as
a single integral equation, namely
\begin{eqnarray}\label{Pi^k}
  \Pi_t^n & = & \Pi_0^n+ \int_0^{t}  \Pi_s^n \Lambda ds+ \sum_{0 < t_k^n \leq t} \left( H_{Y_{t^n_{k}}}[\Pi_{t^n_{k}-}^n]  -  \Pi_{t^n_k-}^n \right),
  \qquad t \geq 0.
     \end{eqnarray}

In what follows we will often use the fact that
 $\Pi^n_t$ is continuously differentiable
on each interval $[t_{k-1}^n,t_{k}^n)$, $k
\geq 0 $ and moreover
 from (\ref{eq-fil-app}) it follows
that there exists a constant $C>0$ (depending only on $\Lambda$) such
that for all $\omega \in \Omega$ and for all $t \notin
\mathcal{P}_n$
\begin{equation}\label{Pn'-bdd}|(\Pi_t^n)'(\omega)| =|\Pi_t^n(\omega) \Lambda| \leq C.
\end{equation}

\noindent {\bf Proof of Theorem \ref{teo-filtro}.}
We consider again the
jumps times $T_j$ ($j \geq 1$) of the process $Y$ and
for every $n$ define $\bar{t}^n_j= \min \{t^n_l \in
\mathcal{P}_n: T_j < t^n_l \}$ and $\underline{t}^n_j= \max
\{t^n_l \in \mathcal{P}_n: t^n_l < T_j\}$ on $\{T_j<\infty\}$;
$\bar{t}^n_j=\underline{t}^n_j=\infty$  on $\{T_j=\infty\}$.
Note that, $\P$-a.s.,
no time $T_j$ belongs to $\mathcal{P}_n$ for $j \geq 1$.
Then, for sufficiently large $n$ (depending on $\omega$),
we have $
\underline{t}^n_j <T_j <\bar{t}^n_j <
\underline{t}^n_{j+1}<T_{j+1}< \bar{t}^n_{j+1}$
on $\{T_j<\infty\}$. We also define
 $T_0=\bar{t}_0^n=\underline{t}^n_0=0$.

For each $j\geq 0$ we consider the following statements:
\begin{itemize}
    \item[$a_j$)] $\Pi^n_{\bar{t}^n_j} \rightarrow \Pi_{T_j}$, $\P$-a.s. on $\{T_j<\infty\}$;
    \item[$b_j$)]  we have, $\P$-a.s.,
    $$
    \begin{array}{l}\dis
    \sup_{t \in [\bar{t}^n_j,T_{j+1})} |\Pi^n_t - \Pi_t| \rightarrow 0
    \mbox { on } \{T_{j+1}<\infty\},\\
    \dis
    \sup_{t \in [\bar{t}^n_j,T)} |\Pi^n_t - \Pi_t| \rightarrow 0
    \mbox { for every } T>0 \mbox{ on } \{T_j<\infty,T_{j+1}=\infty\};\\
    \end{array}
    $$
    \item[$c_j$)]  for all $t \geq 0$ and $i\in I$
    \begin{equation}\label{cj}\Pi_t (i) 1_{T_j\leq t <T_{j+1}}
    =\P(X_t=i|\mathcal{Y}_t^0)1_{T_{j}\leq t <T_{j+1}} \quad \P\mbox{-a.s. }
\end{equation}
\end{itemize}
Note that $a_0$ trivially holds, since
$\Pi_0^n=H_{Y_0}[\mu]=\Pi_0$.
We will prove that the following
implications hold  for all $j \geq 0$:
\begin{eqnarray}
&&a_j \Rightarrow b_j \label{impl1}\\
&&b_j \Rightarrow c_j\label{impl2}\\
&& b_j +c_j\Rightarrow a_{j+1} \label{impl3}
 \end{eqnarray}
 By induction on $j$ it follows in particular that  $c_j$
holds for all $j \geq 0$. This implies that   for all $i \in I$
and $  t \geq 0 $ we have
$$ \Pi_t(i)= \P(X_t=i|\mathcal{Y}_t^0) \qquad
\P-\mbox{a.s. }$$
and concludes the proof of the theorem.

Now it remains to prove (\ref{impl1})-(\ref{impl3}).

\medskip

 \textbf{Proof of (\ref{impl1}).}
 Suppose that $a_j$ is verified for some $j \geq0$. First fix $\omega$
 such that $T_{j+1}(\omega)<\infty$.
 On the interval $[T_j,T_{j+1})$ the process $Y$ is constant.
 For short, we denote $a:=Y_t$ for all $t\in [T_j,T_{j+1})$.
It follows from (\ref{Pn3}) that $\Pi_{t^n_k}^n$ is supported in
$h^{-1}(a)$ for ${t^n_k}\in [\bar{t}^n_j,
\bar{t}^n_{j+1})$.

First take $i \notin h^{-1}(a)$. For $t \in [\bar{t}^n_j,
\bar{t}^n_{j+1})$ there exist two points $t^n_k, t^n_{k+1} \in
\mathcal{P}_n$ such that $\bar{t}^n_j \leq t^n_k \leq t <
t^n_{k+1} \leq \bar{t}^n_{j+1}$. By (\ref{eq-fil-app}) we have
\begin{eqnarray} \Pi_t^n(i) & = &  \Pi_{t^n_k}^n(i)
+\int_{t^n_k}^{t}  (\Pi_s^n \Lambda)(i) \,ds.
 \end{eqnarray}
But $\Pi_{t_k^n}^n(i)=0$ for each $ t_k^n \in [\bar{t}^n_{j}
,T_{j+1})$; hence $\Pi_{t}^n(i)= \int_{t_{k}^n}^{t} (\Pi_s^n
\Lambda)(i) ds$ and from (\ref{Pn'-bdd}) it follows that
$|\Pi_t^n(i)| \leq C(t-t_{k}^n) \leq C 2^{-n}$. We conclude that
 $\sup_{t \in [\bar{t}^n_j,\bar{t}^n_{j+1})} \Pi^n_t(i)\to 0$
as $n\to \infty$ for $i \notin h^{-1}(a)$, or equivalently
\begin{equation}\label{svanisce}
\sup_{t \in [\bar{t}^n_j,\bar{t}^n_{j+1})}
    |  1_{I\backslash h^{-1}(a)}* \Pi_t^n| \to 0,
\end{equation}
which also implies
\begin{equation}\label{svaniscedue}
\sup_{t \in [\bar{t}^n_j,T_{j+1})}
    |  1_{I\backslash h^{-1}(a)}* (\Pi_t^n-\Pi_t)| \to 0,
\end{equation}
since $\Pi_t(i)=0$ for
$i \notin h^{-1}(a)$ and
$t \in [T_j,T_{j+1})\subset [\bar{t}^n_j,T_{j+1})$.

To study $\Pi_{t}^n(i)$
in the case $i \in h^{-1}(a)$, we define for $t \in [\bar{t}^n_j, \bar{t}^n_{j+1})$
$$m^n(t)= \frac{1}{\sum_{j \in h^{-1}(a)} \Pi_t^n(j)}= \frac{1}{ \Pi_t^n 1_{ h^{-1}(a)}} \,.$$
Since $\sum_{i \in I}\Pi_t^n(i)
\equiv 1$ for all $t\geq 0$, by (\ref{svanisce})
we have
$$\sup_{t \in [\bar{t}^n_j,\bar{t}^n_{j+1})} \left|
\sum_{j \in h^{-1}(a)} \Pi_t^n(j)
-1\right| \rightarrow 0,
$$
so that  $m^n(t)$ is well
defined and
\begin{equation}\label{mn}\sup_{t \in [\bar{t}^n_j,\bar{t}^n_{j+1})} |m^n(t)-1| \rightarrow 0.
\end{equation}
We also note that
$m^n$ is cadlag on  $[\bar{t}^n_j,
\bar{t}^n_{j+1})$, and
on each interval $[{t}^n_k,
{t}^n_{k+1})$ contained in $[\bar{t}^n_j,
\bar{t}^n_{j+1})$ it is continuously differentiable  and satifies
$m^n(t_k^n)=1$.

On the interval $[\bar{t}^n_j, \bar{t}^n_{j+1})$, if $i \in h^{-1}(a)$,
$\Pi_{t}^n(i)$ satisfies
\begin{eqnarray*}
  \Pi_t^n(i) & = &  \Pi_{\bar{t}^n_j}^n(i) + \int_{\bar{t}^n_j}^{t}  (\Pi_s^n \Lambda)(i) ds +\sum_{\bar{t}^n_j < t_k^n \leq t} \left( H_a[\Pi_{t^n_{k}-}^n](i)  -  \Pi_{t^n_k-}^n(i) \right)\\
  & = & \Pi_{\bar{t}^n_j}^n(i) + \int_{\bar{t}^n_j}^{t}  (\Pi_s^n \Lambda)(i) ds +
  \sum_{\bar{t}^n_j < t_k^n \leq t} \left(  \frac{\Pi_{t^n_{k}-}^n(i)}{\Pi_{t^n_{k}-}^n 1_{h^{-1}(Y_a)}}  -  \Pi_{t^n_k-}^n(i) \right)\\
  & = & \Pi_{\bar{t}^n_j}^n(i) + \int_{\bar{t}^n_j}^{t}  (\Pi_s^n \Lambda)(i) ds
  +\sum_{\bar{t}^n_j < t_k^n \leq t}  \Pi_{t_k^n-}^n(i)(m^n(t_{k}^n-)-m^n(t_{k-1}^n)) =\\
   &=& \Pi_{\bar{t}^n_j}^n(i) + \int_{\bar{t}^n_j}^{t}  (\Pi_s^n \Lambda)(i) ds +\int_{\bar{t}^n_j}^t \sum_{{\bar{t}^n_j} < t_k^n \leq t}
    \Pi_{t_k^n-}^n(i) 1_{[t_{k-1}^n,t_k^n)}(s)(m^n)'(s) ds. \\
\end{eqnarray*}
In the third equality we used the fact that $m^n(t_k^n)=1$ for
$t_k^n \in [\bar{t}^n_j, \bar{t}^n_{j+1})$. Adding and subtracting
$\int_{\bar{t}^n_j}^t\Pi_s^n(i)(m^n)'(s)ds$ we get
\begin{equation}\label{Ak}\Pi_t^n(i)= \Pi_{\bar{t}^n_j}^n(i)
+ \int_{\bar{t}^n_j}^{t}  (\Pi_s^n \Lambda)(i) ds +
\int_{\bar{t}^n_j}^t\Pi_s^n(i)(m^n)'(s)ds + A_t^n(i)\end{equation}
 where
$$A_t^n(i)=\int_{\bar{t}^n_j}^t \left[\sum_{\bar{t}^n_j < t_k^n \leq t}
\Pi^n_{t_k^n-}(i)  1_{[t_{k-1}^n,t_k^n)}(s)
-\Pi_s^n(i)\right](m^n)'(s)\,ds.
$$
Thanks to
(\ref{Pn'-bdd}), which also implies that $|(m^n)'(s)| \leq C$ for
a.e. $s$, we obtain
\begin{equation}\label{A}\sup_{t\in [\bar{t}^n_j , \bar{t}^n_{j+1})} |A_t^n(i)| \rightarrow 0
 \mbox{ as } n \rightarrow \infty.
\end{equation}
Moreover we can compute
$$(m^n)'(t)= -\frac{\sum_{j \in h^{-1}(a)}
\Pi_t^n(j)'}{(\sum_{j \in h^{-1}(a)} \Pi_t^n(j))^2}= -\sum_{j
\in h^{-1}(a)} \Pi_t^n(j)'+B^n(t),$$ where
$$B^n(t)=\sum_{j \in h^{-1}(a)} \Pi_t^n(j)'\left(1-
\frac{1} {(\sum_{j \in h^{-1}(a)} \Pi_t^n(j))^2}\right)
=
\sum_{j \in h^{-1}(a)} \Pi_t^n(j)'\left(1-
m^n(t)^2\right)
$$
and we have
 \begin{equation}\label{B}\sup_{t\in [\bar{t}^n_j,
\bar{t}^n_{j+1})} |B^n(t)| \rightarrow 0 \mbox{  as } n
\rightarrow \infty
\end{equation}
by (\ref{Pn'-bdd}) and (\ref{mn}).
Now (\ref{Ak}) becomes
$$
\begin{array}{lll}
\Pi_t^n(i)& = &\dis \Pi_{\bar{t}^n_j}^n(i) + \int_{\bar{t}^n_j}^{t}  (\Pi_s^n \Lambda)(i) ds
- \int_{\bar{t}^n_j}^t\Pi_s^n(i)\sum_{j \in h^{-1}(a)} \Pi_s^n(j)'ds +\\
& & \dis+
\int_{\bar{t}^n_j}^t\Pi_s^n(i) B^n(s)ds + A_t^n(i)\\
& = & \dis\Pi_{\bar{t}^n_j}^n(i) + \int_{\bar{t}^n_j}^{t}  (\Pi_s^n
\Lambda)(i) ds -\int_{\bar{t}^n_j}^t\Pi_s^n(i)\sum_{j \in
h^{-1}(a)} \Pi_s^n(j)'ds + C^n_t(i) + A_t^n(i)
\end{array}
$$
where we have defined
\begin{equation}\label{defcn}
    C^n_t(i)= \int_{\bar{t}^n_j}^t\Pi_s^n(i) B^n(s)ds.
\end{equation}
By (\ref{eq-fil-app}), for $ s \in [\bar{t}^n_j,t)$, $s\notin \calp_n$, we have
 $\sum_{j \in h^{-1}(a)}
\Pi_s^n(j)'= \sum_{j \in h^{-1}(a)} (\Pi_s^n \Lambda)(j)=
\Pi_s^n \Lambda 1_{h^{-1}(a)}$ and we get
$$\Pi_t^n(i) =  \Pi_{\bar{t}^n_j}^n(i) + \int_{\bar{t}^n_j}^{t}  (\Pi_s^n \Lambda)(i) ds
- \int_{\bar{t}^n_j}^t\Pi_s^n(i)(\Pi_s^n \Lambda 1_{h^{-1}(a)} )ds
+ C_t^n(i) + A_t^n(i).
$$
We denote with $A^n_t$ and $C^n_t$ the $N$-dimensional
vectors with components $A^n_t(i)$ and $C^n_t(i)$, respectively,  if
$i \in h^{-1}(a)$, and $0$ otherwise.
Introducing  the vector field $G_a$ defined  by
$$
G_a(\nu)=1_{h^{-1}(a)}*[ (\nu\Lambda) - (\nu\Lambda
1_{h^{-1}(a)})\nu],\qquad \nu\in \Delta
$$
we finally arrive at
$$
1_{h^{-1}(a)}*\Pi_t^n = 1_{h^{-1}(a)}* \Pi_{\bar{t}^n_j}^n +
\int_{\bar{t}^n_j}^{t} G_a( \Pi_s^n)\, ds
+ C_t^n + A_t^n,\qquad t\in  [\bar{t}^n_j, T_{j+1}).
$$

Now we are able to obtain an estimate on the difference
$1_{h^{-1}(a)}*[\Pi_t-\Pi_t^n]$ for $t \in [\bar{t}^n_j,
T_{j+1})$. We recall that $\Pi_t \in \Delta_a$ for all $t \in [T_j,
T_{j+1})$, so that
(\ref{eqdiffintervallo}) implies
$$
1_{h^{-1}(a)}*\Pi_t= 1_{h^{-1}(a)}* \Pi_{T_j} +
\int_{T_j}^{t} G_a( \Pi_s)\, ds, \qquad t\in  [T_j, T_{j+1}).
$$
So, for $t \in [\bar{t}^n_j,T_{j+1})$,
\begin{eqnarray*}
|1_{h^{-1}}(a)*(\Pi_t- \Pi_t^n)| & \leq  & |1_{h^{-1}}(a)*(\Pi_{T_j}- \Pi_{\bar{t}^n_j}^n)|+
\int_{T_j}^{\bar{t}^n_j}|G_a( \Pi_s)|\, ds +
|C_t^n| + |A_t^n|\\
 &  &+\int_{\bar{t}^n_j}^t | G_a( \Pi_s)-G_a( \Pi_s^n)|\,ds
 \end{eqnarray*}
Note that $G_a$
is bounded and globally Lipschitz on $\Delta$. We denote $K$ some bound and
by $L$ its
Lipschitz constant ($K$ and $L$ depend on $\omega$, since
$a=Y_{T_j}(\omega)$).
Then we obtain, noting that ${\bar{t}^n_j}-T_j\le 2^{-n}$
$$
\begin{array}{l}
\dis
|1_{h^{-1}}(a)*(\Pi_{T_j}- \Pi_{\bar{t}^n_j}^n)|+
\int_{T_j}^{\bar{t}^n_j}|G_a( \Pi_s)|\, ds +
|C_t^n| + |A_t^n|
\\
\qquad \dis
 \leq
|1_{h^{-1}}(a)*(\Pi_{T_j}- \Pi_{\bar{t}^n_j}^n)|+
K 2^{-n}+
|C_t^n| + |A_t^n|
  \end{array}
  $$
and
\begin{eqnarray*}
|1_{h^{-1}}(a)*(\Pi_t- \Pi_t^n)| & \leq  &
|1_{h^{-1}}(a)*(\Pi_{T_j}- \Pi_{\bar{t}^n_j}^n)|+
K 2^{-n}+
|C_t^n| + |A_t^n|
\\
 &  &+L\int_{\bar{t}^n_j}^t |  \Pi_s-\Pi_s^n|\,ds.
 \end{eqnarray*}
Since $1_{h^{-1}(a)}*\Pi_s=\Pi_s$ we have
$
|  \Pi_s-\Pi_s^n|\le
|1_{h^{-1}(a)}*(\Pi_s- \Pi_s^n)| +
|1_{I\backslash h^{-1}(a)}*\Pi_s|
$
and we obtain
\begin{eqnarray*}
|1_{h^{-1}(a)}*(\Pi_t- \Pi_t^n)| & \!\!\leq \!\! &
|1_{h^{-1}(a)}*(\Pi_{T_j}- \Pi_{\bar{t}^n_j}^n)|+
K 2^{-n}+
|C_t^n| + |A_t^n|+
L\int_{\bar{t}^n_j}^t | 1_{I\backslash h^{-1}(a)}* \Pi_s^n|\,ds
\\
 &  &+L\int_{\bar{t}^n_j}^t | 1_{h^{-1}(a)}*(\Pi_s- \Pi_s^n)|\,ds
\\
 &  &=:D_t^n
 +L\int_{\bar{t}^n_j}^t | 1_{h^{-1}(a)}*(\Pi_s- \Pi_s^n)|\,ds.
 \end{eqnarray*}

From the induction assumption $a_j$ and from
(\ref{svanisce}), (\ref{defcn}), (\ref{A}),  (\ref{B})
it follows that
 $\sup_{t\in [\bar{t}^n_j, T_{j+1})}|D^n_t|\to 0$. From the
Gronwall lemma we  conclude that
 $\sup_{t\in [\bar{t}^n_j, T_{j+1})}|1_{h^{-1}}(a)*(\Pi_t- \Pi_t^n)|\to 0$.
 Recalling (\ref{svaniscedue}), this
proves  $b_j$ in the case $T_{j+1}(\omega)<\infty$.

 The case $ T_j(\omega)<\infty,T_{j+1}(\omega)=\infty$ can be proved
 by the same arguments, replacing the interval
 $[\bar{t}^n_j, T_{j+1})$ by $[\bar{t}^n_j, T)$ in the previous passages.

\medskip

\textbf{Proof of (\ref{impl2}).} Assume that $b_j$ holds for some $j\ge 0$.
For every $i\in I$ and $t\ge 0$, by
lemma \ref{proprietaPn}
 we have
$$\Pi^n_t(i)\, 1_{T_j < t
<T_{j+1}}=\P(X_t=i|\mathcal{Y}^n_t)\,1_{T_{j} < t <T_{j+1}},
\qquad
\P \mbox{-a.s.}
$$
Let $n\to\infty$.  Since $\bar{t}^n_j\to T_j$, for large $n$ we have
$\bar{t}^n_j\le t$ and $b_j$ implies that $\Pi^n_t(i)\to \Pi_t(i)$.
Recalling (\ref{convfiltrovero}) we conclude that
$$\Pi_t(i)\, 1_{T_j < t
<T_{j+1}}=\P(X_t=i|\mathcal{Y}^0_t)\,1_{T_{j} < t <T_{j+1}},
\qquad
\P \mbox{-a.s.}
$$
Since $\P(T_j=t)=0$ this shows that $c_j$ holds.

\medskip

\textbf{Proof of (\ref{impl3}).}
We suppose that $b_j$ and $c_j$ hold for some $j \geq 0$.
We write
$\Pi^n(i,t)$ and $\Pi^n(t)$ instead of
 $\Pi_t^n(i)$ and $\Pi_t^n$ for better readability.
We recall that, by (\ref{eqdiffintervallo}), for
 $ T_{j}\le
t <T_{j+1}$, we have $\Pi_t(i)=0$ if $i\notin h^{-1}(Y_{T_{j}})$ and
$$
     \Pi_t^{'}(i)=
 (\Pi_t\Lambda)(i) -
(\Pi_t\Lambda 1_{h^{-1}(Y_{T_{j}})})\Pi_t(i)
$$
if $i\in h^{-1}(Y_{T_{j}})$. Since we also assume that
 $c_j)$ holds, Lemma
\ref{P-Tj} ensures that
\begin{equation}\label{nonsiannulladen}
\Pi(T_{j+1}-)\Lambda 1_{h^{-1}
(Y_{T_{j+1}})}>0 \mbox{ on  } \{T_{j+1}<\infty\}.
\end{equation}

Now we fix $\omega$ such that $T_{j+1}<\infty$ and we wish to prove that
$\lim_{n \rightarrow
\infty} \Pi^n(\bar{t}^n_{j+1})=
 \Pi_{T_{j+1}}$. Note that $Y_{T_{j+1}}=Y_{\bar{t}^n_{j+1}}$ so that
if $i \notin  h^{-1} (Y_{T_{j+1}}) = h^{-1} (Y_{\bar{t}^n_{j+1}})$, then
$\Pi^n(i,\bar{t}^n_{j+1})=0= \Pi_{T_{j+1}}(i)$. Let now $i \in   h^{-1}
(Y_{T_{j+1}})$; then, since by (\ref{Pn3}) we have
$\Pi^n(\bar{t}^n_{j+1})=
H_{Y_{\bar{t}^n_{j+1}}}[\Pi^{
n}(\bar{t}^n_{j+1}-)]$, we therefore have
\begin{equation}\label{dafarelimite}
    \Pi^n(i,\bar{t}^n_{j+1})
=
\frac{\Pi^{
n}(i,\bar{t}^n_{j+1}-)}{\sum_{k \in h^{-1} (Y_{\bar{t}^n_{j+1}})}
\Pi^{ n}(k,\bar{t}^n_{j+1}-) }
=\frac{\Pi^{
n}(i,\bar{t}^n_{j+1}-)}{\sum_{k \in h^{-1} (Y_{T_{j+1}})}
\Pi^{ n}(k,\bar{t}^n_{j+1}-) },
\end{equation}
provided the denominator is not zero, and we will check that this is indeed the case.
Thanks to  (\ref{eq-fil-app}) we can compute the right-hand side by a
Taylor's expansion around $\underline{t}^n_{j+1}$:
\begin{eqnarray*}
\Pi^n(\bar{t}^n_{j+1}-)&= &
\Pi^n(\underline{t}^n_{j+1})+
(\bar{t}^n_{j+1}-\underline{t}^n_{j+1})\;(\Pi^{ n})^{'}(\underline{t}^n_{j+1})+
 o(\bar{t}^n_{j+1}-\underline{t}^n_{j+1})
 \\
 &=&
\Pi^n(\underline{t}^n_{j+1})+
(\bar{t}^n_{j+1}-\underline{t}^n_{j+1})\;\Pi^{ n}(\underline{t}^n_{j+1})\Lambda+
 o(\bar{t}^n_{j+1}-\underline{t}^n_{j+1}).
\end{eqnarray*}
Noting that $\bar{t}^n_{j+1}-\underline{t}^n_{j+1}=
2^{-n}$ and that
$\Pi^n(i,\underline{t}^n_{j+1})=\Pi^n(k,\underline{t}^n_{j+1})=0$ for $k \in h^{-1} (Y_{T_{j+1}})$
we obtain
$$\Pi^n(i,\bar{t}^n_{j+1})=
\frac{2^{-n}\;
[\Pi^n(\underline{t}^n_{j+1}) \Lambda
](i)+o(2^{-n})}{2^{-n}\;
 \sum_{k \in h^{-1} (Y_{T_{j+1}})}[\Pi^n(\underline{t}^n_{j+1})\Lambda ](k) +o(2^{-n})}.
 $$
Since $\underline{t}^n_{j+1}$ tends to $T_{j+1}-$, and since  $b_j$
shows that $\Pi^n(t)$ converges to $\Pi_t$ uniformly in a left
neighborhood of $T_{j+1}$, we have $\Pi^n(\underline{t}^n_{j+1})\to
\Pi(T_{j+1}-)$
and we finally obtain
\begin{equation}\label{PnTj} \lim_{n \rightarrow \infty}\Pi^n(i,\bar{t}^n_{j+1})
=\frac{  [\Pi(T_{j+1}-) \Lambda ](i)}{
 \sum_{k \in h^{-1} (Y_{T_{j+1}})}[\Pi(T_{j+1}-)\Lambda ](k) }= \frac{  [\Pi(T_{j+1}-) \Lambda ](i)}{
\Pi(T_{j+1}-)\Lambda  1_{h^{-1} (Y_{T_{j+1}})}}.
\end{equation}
We note that the
right-hand side of (\ref{PnTj})
 is well defined by (\ref{nonsiannulladen}). For the same
 reason the denominator in (\ref{dafarelimite}) is strictly positive
 for large $n$, which justifies previous passages.
 Finally, since
$\Pi_{T_{j+1}}=H_{Y_{T_{j+1}}}[\Pi(T_{j+1}-) \Lambda]$ by
(\ref{eq-pifiltro-equiv}),  (\ref{PnTj}) also shows that
$\lim_{n \rightarrow \infty}\Pi^n(i,\bar{t}^n_{j+1})=
\Pi_{T_{j+1}}(i)$ and concludes the proof that
$a_{j+1}$ holds. \qed

\section{Basic properties of the filtering process}
\label{basicprop}

\subsection{Canonical set-up}\label{canonical}

In the rest of this paper it is convenient to
assume that the unobserved process is defined in a canonical set-up
as follows.

\begin{enumerate}
\item Let $\Omega$  be the set of cadlag functions $\omega:\R_+\to I$,
i.e. the set of right-continuous functions having finite left limits on $(0,\infty)$.
We denote $X_t(\omega)=\omega(t)$ for $\omega\in \Omega$ and $t\ge 0$,
and we introduce the $\sigma$-algebras
$$
\calf^0_t=\sigma(X_s\,:\,s\in[0,t]),
 \qquad
\calf^0=\sigma(X_s\,:\,s\ge0).
$$
 $(\calf^0_t)_{t\ge0}$ is thus the natural filtration of $(X_t)_{t\ge0}$.

\item Let $\Delta$  denote the set of probability measures on $I$, identified
 with the canonical simplex of $\R^N$, where $N$ is
the cardinality of $I$.

\item For every $\mu\in\Delta$ we denote by $P_\mu$ the unique probability measure on
$(\Omega,\calf^0)$ that makes $(X_t)$ a Markov process on $I$ with generator
$\Lambda$ and initial distribution $\mu$, i.e. such that
for every $t,s\ge0$ and every real function $f$ on $I$ we have
$$
E_\mu [f(X_{t+s})|\calf^0_t]=(e^{s\Lambda}f)(X_t),\qquad P_\mu-a.s.
$$
and $P_\mu(X_0=i)=\mu(i)$, $i\in I$. Here $E_\mu$ denotes of course the
expectation with respect to $P_\mu$.

If $\mu$ is concentrated at some $i\in I$ we write $P_i$ instead of $P_\mu$.

\item We still define the observation process $(Y_t)_{t\ge0}$
and  its
natural filtration  $(\caly^0_t)_{t\ge0}$ by
$$Y_t=h(X_t),\quad
\caly^0_t=\sigma(Y_s\,:\,s\in[0,t]),\qquad t\ge 0.
$$

\end{enumerate}

\begin{remark}\begin{em}
We could also define the space $\Omega$ as the set of all functions
$\omega:\R_+\to I$ which are
piecewise-constant, right-continuous and
with a finite number of jumps in every bounded interval, i.e. of the form
$\omega(t)=\sum_{k=0}^\infty a_k 1_{[t_k,t_{k+1})}(t)$ for $a_k\in\R$
 and for $0=t_0<t_1<t_2<\ldots$ where the sequence $(t_k)$ has no
 cluster point in $[0,\infty)$. The other definitions remain unchanged, and all
 the subsequent results still hold.
\end{em}
  \end{remark}

Let $(T_n)_{n\ge 1}$ denote the  sequence of  jumps times of $(Y_t)$,
with the convention that $T_n=\infty$ for all $n$ if no jump occurs
and $T_{n}<T_{n+1}=T_{n+2}=\ldots=\infty$ if precisely $n$ jumps occur.

For every  $\omega\in\Omega$ we consider the corresponding trajectory $Y_t(\omega)$
and jump times
 $T_n(\omega)$
and we define the filtering process $\Pi^\mu_t(\omega)$ as the solution of
\begin{equation}\label{eqdefpifiltro}
\begin{array}{lll}
\dis
\Pi^\mu_t(\omega)&=&\dis
H_{Y_0(\omega)}[\mu]+\int_0^t\left\{
1_{h^{-1}(Y_s(\omega))}\ast (\Pi^\mu_s(\omega)\Lambda)
- (\Pi^\mu_s(\omega)\Lambda 1_{h^{-1}(Y_s(\omega))})\Pi^\mu_s(\omega)\right\}ds
\\&&\dis
+\sum_{0<T_n(\omega)\le t}\left\{
H_{Y_{T_n}(\omega)}[\Pi^\mu_{T_n-}(\omega)\Lambda]-\Pi^\mu_{T_n-}(\omega)\right\},
\qquad t\ge0,
\end{array}
\end{equation}
where, as usual, $\Pi^\mu_{T_n-}(\omega)
=\lim_{t\to {T_n}(\omega), t<{T_n}(\omega)}\Pi^\mu_t(\omega)$
is defined on $\{\omega\in\Omega\,:\,{T_n}(\omega)<\infty\}$.

 Equation (\ref{eqdefpifiltro})
is the same as (\ref{eq-pifiltro}) and, as explained before,
uniquely
determines  a $(\caly^0_t)$-adapted, cadlag process $(\Pi^\mu_t)$
taking values in the effective simplex $\Delta_e=\cup_{a\in O}\Delta_a$.
By theorem  \ref{teo-filtro}, for every $\mu\in\Delta$, $t\ge0$ and $i\in I$
 we have $\Pi^\mu_t(i)=P_\mu(X_t=i|\caly_t^0)$,
$P_\mu$-a.s. and,  consequently,
$\Pi^\mu_t(i)=P_\mu(X_t=i|\caly_t^\mu)$,
$P_\mu$-a.s.

In what follows the process $(\Pi^\mu_t)$ will also be considered under a different
probability $P_\rho$ with $\rho\in\Delta$, $\rho\neq \mu$. In this case
the equality $\Pi^\mu_t(i)=P_\rho(X_t=i|\caly_t^0)$ is generally false.

\begin{remark}\begin{em}\label{piepinu}
We note that
the process $(\Pi^\mu_t)$ depends on $\mu$ through its initial value
$\Pi^\mu_0= H_{Y_0}[\mu]$. In general  $\Pi^\mu_0$  is
 a random variable and
 $\Pi^\mu_0\neq \mu$.
However if the unobserved process $(X_t)$ has initial distribution $\nu$ in the
effective simplex $\Delta_e$ then
 the filtering process $(\Pi^\nu_t)$ starts at $\nu$, $P_\nu$-a.s.
 Indeed, if $\nu\in\Delta_a$, then
 $Y_0=a$ $P_\nu$-a.s. and since
 $H_a[\nu]=\nu$ it follows that $\Pi^\nu_0=\nu$ $P_\nu$-a.s.
\end{em}
  \end{remark}

\subsection{Prediction}

Preliminary to further properties of the filtering process we need to prove the
following result which has an intrinsic interest, since it solves
the so called prediction problem: at any time $t$
it allows to compute the distribution of the unobserved process
after time $t$ conditional on the available observation up to $t$.

\begin{proposition}\label{prediction}
For every $\mu\in\Delta$, $t\ge0$ and $\Gamma\in \calf^0$ we have
$$
P_\mu(X_{t+\cdot}\in\Gamma|\caly_t^0)=P_{\Pi^\mu_t}(\Gamma),
\qquad P_\mu-a.s.
$$
\end{proposition}

In this statement $P_{\Pi^\mu_t}(\Gamma)$ denotes
$\sum_{i\in I}P_{i}(\Gamma)\,\Pi^\mu_t(i)$.
Roughly, the proposition states that for every present time
$t$, and conditional on
the past observation of $Y$ up to $t$,
the probability that the future trajectories of $X$ will belong
to some set $\Gamma$ is
best predicted by $P_{\Pi^\mu_t}(\Gamma)$, i.e. one computes
$P_{\mu}(\Gamma)$ and replaces $\mu$ by ${\Pi^\mu_t}$.

 \noindent {\bf Proof.} Noting that $\caly_t^0\subset \calf_t^0$, by the
 Markov property of $X$ and the fact that $\Pi^\mu$ is the filtering process
 we have
 $$
 P_\mu(X_{t+\cdot}\in\Gamma|\caly_t^0)
 =
E_\mu[ P_\mu(X_{t+\cdot}\in\Gamma|\calf_t^0)|\caly_t^0]
=E_\mu[ P_{X_t}(\Gamma)|\caly_t^0]
=\sum_{i\in I}P_{i}(\Gamma)\,\Pi^\mu_t(i). \qquad\qed
 $$

 \subsection{The Markov property of  the filtering process}

For $t\ge0$, $\nu\in\Delta_e$  and $A\in \calb(\Delta_e)$ (the Borel $\sigma$-algebra
of $\Delta_e$) we define
$$
R_t(\nu,A)=P_\nu(\Pi_t^\nu\in A).
$$
It can be easily verified that $R_t$ is a Markov kernel on $(\Delta_e,\calb(\Delta_e))$.
The following proposition asserts
the Markov property of the filtering process $(\Pi^\mu_t)$, corresponding
to arbitrary fixed initial distribution $\mu\in\Delta$ of the unobserved process $(X_t)$.

\begin{proposition}\label{filtromarkoviano}
$(R_t)$ is a Markov transition function on $(\Delta_e, \calb(\Delta_e))$ and
for  $\mu\in\Delta$, $t,s\ge0$ and $A\in \calb(\Delta_e)$ we have
$$
P_\mu(\Pi^\mu_{t+s}\in A|\caly_t^0)=R_s(\Pi^\mu_t,A),
\qquad P_\mu-a.s.
$$
In other words, for every $\mu\in\Delta$, in the probability space $(\Omega, \calf^0,P_\mu)$ the process
$(\Pi^\mu_t)$ is a Markov process with respect to $(\caly_t^0)$,
taking values in $\Delta_e$ and having transition function $(R_t)$.
\end{proposition}

 \noindent {\bf Proof.}
We first introduce a family of stochastic processes
$(\Pi_t(\nu))$ parametrized by $\nu\in\Delta_e$. For every $\omega\in\Omega$
we define  $\Pi_t(\omega,\nu)$ as the solution of
\begin{equation}\label{eqdefpifiltronu}
\begin{array}{lll}
\Pi_t(\omega,\nu)&=&\dis
\nu+\int_0^t\left\{
1_{h^{-1}(Y_s(\omega))}\ast (\Pi_s(\omega,\nu)\Lambda)
- (\Pi_s(\omega,\nu)\Lambda 1_{h^{-1}(Y_s(\omega))})\Pi_s(\omega,\nu)\right\}ds
\\&&\dis
+\sum_{0<T_n(\omega)\le t}\left\{
H_{Y_{T_n}(\omega)}[\Pi_{T_n-}(\omega,\nu)\Lambda]-\Pi_{T_n-}(\omega,\nu)\right\},
\qquad t\ge0.
\end{array}
\end{equation}
Since $\nu$ belongs to $\Delta_e$,
 $(\Pi_t(\nu))$ takes values in  $\Delta_e$
and it is a $(\caly^0_t)$-adapted, cadlag process.
Moreover $(\omega,t,\nu)\to \Pi_t(\omega,\nu)$ is
measurable with respect to $\calf^0\times \calb(\R_+)\times \calb(\Delta_e)$.

Note that  the pathwise evolution of $(\Pi_t(\nu))$ is described by the same
differential equation as for $(\Pi^\mu_t)$, but these two processes
differ in general because of the initial condition:
$(\Pi_t(\nu))$ starts at $\nu$, whereas $\Pi^\mu_0(\omega)= H_{Y_0(\omega)}[\mu]$.
Clearly, we have $\Pi^\mu_t(\omega)=\Pi_{t}(\omega,\Pi^\mu_0(\omega ))$.
If $\mu=\nu\in\Delta_e$ then, as noted in remark \ref{piepinu}, we have $\Pi^\nu_0=\nu$
and consequently $(\Pi^\nu_t)$ and $(\Pi_t(\nu))$ are the same process.

Let us define the translation operators $\theta_t:\Omega\to\Omega$ by
$(\theta_t\omega)(s)=\omega(t+s)$ for $t,s\ge 0$. By the uniqueness of
the solution of equation
(\ref{eqdefpifiltronu}) we have,
for all $t,s\ge 0$ and $\omega\in\Omega$,
$$
\Pi_{t+s}(\omega,\nu)=\Pi_{s}(\theta_t\omega, \Pi_{t}(\omega,\nu))
$$
and replacing $\nu$ by $\Pi^\mu_0(\omega)$ we obtain
$$\Pi_{t+s}^\mu(\omega)=\Pi_{t+s}(\omega,\Pi^\mu_0(\omega))
= \Pi_{s}(\theta_t\omega, \Pi_{t}(\omega,\Pi^\mu_0(\omega) ))
=\Pi_{s}(\theta_t\omega, \Pi_{t}^\mu(\omega )).
$$
Fixing $t,s\ge 0$ and noting that  $\Pi_{t}^\mu$ is $\caly_t^0$-measurable
it follows that
\begin{equation}\label{pimarkov}
    P_\mu\left(\Pi^\mu_{t+s}\in A|\caly_t^0\right)=
P_\mu\left(\Pi_{s}(\theta_t(\cdot), \Pi_{t}^\mu ))\in A|\caly_t^0\right)=
g(\Pi_{t}^\mu),\qquad P_\mu-a.s.,
\end{equation}
where  we define
$$
g(\rho):=
P_\mu\left(\Pi_{s}(\theta_t(\cdot), \rho))\in A|\caly_t^0\right), \qquad \rho \in\Delta_e.
$$
Note that the event $\{\Pi_{s}(\theta_t(\cdot), \rho))\in A\}$ can be written in the form
$\{X_{t+\cdot}\in\Gamma_\rho\}$ where
$\Gamma_\rho:= \{\omega\in\Omega\,:\, \Pi_{s}(\omega, \rho))\in A\}$.
So it follows from proposition \ref{prediction} that
$$
g(\rho)= P_{\mu}\left(X_{t+\cdot}\in\Gamma_\rho|\caly_t^0 \right)
=P_{\Pi^\mu_t}\left(\Gamma_\rho \right)
=P_{\Pi^\mu_t}\left(\Pi_{s}(\rho)\in A \right).
$$
Replacing in (\ref{pimarkov}) we obtain the
required equality
$P_\mu\left(\Pi^\mu_{t+s}\in A|\caly_t^0\right)=
g(\Pi_{t}^\mu)= R_s(\Pi_{t}^\mu,A)$. The Chapman-Kolmogorov equation
for $(R_t)$ now follows from this equality and the fact noted earlier that
$\Pi^\nu_t=\Pi_t(\nu)$ for $\nu\in\Delta_e$.
  \qed

\begin{remark}\label{instabfilter} (Filter instability).
{\em
An important issue, both for theoretical and computational viewpoint, is
the stability of the filtering process. Stability is often formulated in terms of
appropriate
ergodic properties of the Markov filtering process. It is an interesting fact
that stability essentially fails to hold
in the case of noise-free observation under consideration.
This is the content of section 3 of \cite{BaChiLip}, where
the authors consider
 the Markov chain
with rate
transition matrix
 given by (\ref{lambdacontroes}) and the observation process corresponding
 to the function $h$ defined in
(\ref{hcontroes}). They show that the filtering process has infinitely many
invariant measures, and that the solutions of the filtering equations corresponding
to different initializations do not converge to one another: see subsection 3.2
in \cite{BaChiLip} for a detailed discussion.
}
\end{remark}

\section{The filtering process and the distribution of the observation process}
\label{sectionobservation}

In this section we show that the law of the observation process can be described
by means of the filtering process through explicit formulae.
 These results will be applied in section \ref{sectionpdp}, but they
also have some intrinsic interest: we present an  application in subsection
\ref{tempouscita}.

We will use the canonical set-up introduced in subsection
\ref{canonical}, we fix $\mu\in\Delta$ as an
initial distribution of $X$
and we construct
the probability  space $(\Omega,\calf^0,P_\mu)$. In the rest of this
section all stochastic processes
will be considered under $P_\mu$.

Let $(T_n)_{n\ge 1}$ be the  sequence of  jumps times of $Y$,
with the convention that $T_n=\infty$ for all $n$ if no jump occurs
and $T_{n}<T_{n+1}=T_{n+2}=\ldots=\infty$ if precisely $n$ jumps occur. We let
 $T_0=0$. In the following we will consider the sojourn times
  $S_n=T_n-T_{n-1}$  and the positions
  $Y_{T_{n-1}}$,   $Y_{T_{n-1}-}$,
   of $Y$ at jump times and immediately before jump times respectively
($n\ge 1$). These random variables are only defined
   on the event $\{ T_{n-1}<\infty\}$. Nevertheless the $\sigma$-algebras
\begin{equation}\label{defsigma}
    \Sigma_{n-1}=\sigma (Y_0,T_1,Y_{T_1},\ldots, T_{n-1},Y_{T_{n-1}}),\quad
\Sigma^+_{n-1}=\sigma (Y_0,T_1,Y_{T_1},\ldots, T_{n-1},Y_{T_{n-1}},T_n ),
\end{equation}
can be defined in the usual way. In particular
$\Sigma_{0}=\sigma (Y_0)$, $\Sigma^+_{0}=\sigma (Y_0,T_1)$.

Since the trajectories of $Y$ are constant among jump times,
the law of $Y$ is completely
determined by  the finite-dimensional
distributions of the stochastic process
$\{Y_0,T_1,$ $ Y_{T_1},$ $ T_2, $ $Y_{T_2},$ $\ldots\}$. These in turn
can be described by specifying the distribution of $Y_0$,
which is obvious, and the family of conditional
probabilities
$$
P_\mu(S_n>t,T_{n-1}<\infty\,|\,
\Sigma_{n-1}),\quad
P_\mu(Y_{T_n}=b,{T_n<\infty}\,|\,
\Sigma^+_{n-1}),
\qquad t\ge 0,\, b\in O,\, n\ge 1.
$$

In order to present explicit formulae for these probabilities we need to introduce
some notation.
For
$\nu \in\Delta_e$ we define a probability $q(\nu)=(q(\nu,b))_{b\in O}$
on $O$ in the following way. For every $a\in O$ we first fix
an arbitrary probability $q_a=(q_a(b))_{b\in O}$ on $O$ supported in $O\backslash \{a\}$
(we exclude the trivial case where $h$ is constant); the
exact values of $q_a$ are irrelevant.
Next, if $\nu\in\Delta_a$
for some $a\in O$, we define for every $b\in O$
$$
q(\nu ,b)=
\left\{
\begin{array}{ll}
\dis
\frac{\nu \Lambda 1_{h^{-1}(b)}
}{-\nu \Lambda 1_{h^{-1}(a)}
}\,1_{b\neq a}
,
&{\rm if\;} \nu \Lambda 1_{h^{-1}(a)}\neq 0,
\\\dis
q_a(b),
&{\rm if\;} \nu \Lambda 1_{h^{-1}(a)}= 0.
\end{array}
\right.
$$
To check that  $q(\nu)$ is a probability measure
we first note that for $\nu \in\Delta_a$ and $b\neq a$ we have
$$
\nu \Lambda 1_{h^{-1}(b)}=\sum_{i\in h^{-1}(a)}\sum_{j\in h^{-1}(b)}
\nu _i\lambda_{ij}
\ge 0,
$$
since the sums are extended to distinct indices $i,j$ and therefore
$\lambda_{ij}\ge0$.
Moreover for $\nu \in\Delta_a$
$$\nu \Lambda 1_{h^{-1}(a)}+
\sum_{b\in O,b\neq a}\nu \Lambda 1_{h^{-1}(b)}=
\sum_{b\in O}\nu \Lambda 1_{h^{-1}(b)}=\nu \Lambda 1=0,
$$
since $\Lambda 1=0$,
so that $
-\nu \Lambda
1_{h^{-1}(a)}\ge 0$ and $q(\nu ,\cdot)$ is in fact a probability measure.
Note  that if $\nu \in\Delta_a$ and
$\nu \Lambda 1_{h^{-1}(a)}=0$ then also
$\nu \Lambda 1_{h^{-1}(b)}=0$ for all $b\neq a$ and consequently
the equality
\begin{equation}\label{cancq}
-\nu \Lambda 1_{h^{-1}(a)} q(\nu ,b)=
\nu \Lambda 1_{h^{-1}(b)},
\end{equation}
holds for all $\nu \in\Delta_a$ and $b\neq a$.
Finally note  that if $\nu \in\Delta_a$ then $q(\nu,a)=0$.

Let
 $(\Pi^\mu_t)$ denote the
 the filtering process,
 solution
 of (\ref{eqdefpifiltro}).

\begin{theorem}\label{mainyepi}
Let $\mu\in\Delta$. Then
for $t\ge 0$, $b\in O$, $n\ge1$ we have
 \begin{equation}    \label{sncondizy}
P_\mu(S_n>t,T_{n-1}<\infty\,|\,\Sigma_{n-1})=
\exp\left(\int_0^t\phi(s,\Pi^\mu_{T_{n-1}})\Lambda
1_{h^{-1}(Y_{T_{n-1}})}ds\right)\,1_{T_{n-1}<\infty},
 \end{equation}
 \begin{equation}    \label{ytncondizy}
P_\mu(Y_{T_n}=b,T_n<\infty\,|\,\Sigma^+_{n-1} )
=q\left(\Pi^\mu_{T_{n}-},b\right)1_{T_n<\infty}
=q\left(\phi (S_n, \Pi^\mu_{T_{n-1}}),b\right)1_{T_n<\infty}.
 \end{equation}
\end{theorem}

\begin{remark}\begin{em}
Recall that $\Sigma_{n-1}$ and $\Sigma^+_{n-1}$ were defined in
(\ref{defsigma}).
The equalities (\ref{sncondizy})-(\ref{ytncondizy}) may be written
$$
P_\mu(S_n>t\,|\,\Sigma_{n-1})=
\exp\left(\int_0^t\phi(s,\Pi^\mu_{T_{n-1}})\Lambda
1_{h^{-1}(Y_{T_{n-1}})}ds\right)
\;{\rm on\;} \{T_{n-1}<\infty\},
$$
$$
P_\mu(Y_{T_n}=b\,|\,\Sigma^+_{n-1})
=q\left(\Pi^\mu_{T_{n}-},b\right)
=q\left(\phi (S_n, \Pi^\mu_{T_{n-1}}),b\right)
\;{\rm on\;} \{T_n<\infty\}.
$$
The intuitive meaning of (\ref{sncondizy})-(\ref{ytncondizy}) is as follows:
suppose we know
$Y_0,T_1,\Pi^\mu_{T_1},\ldots, T_{n-1},Y_{T_{n-1}}$, or equivalently we know
the trajectory of $Y$ up to the jump time $T_{n-1}$. This also determines
the trajectory of $\Pi^\mu$ up to $T_{n-1}$.
Suppose that $Y_{T_{n-1}}=a$, $\Pi^\mu_{T_{n-1}}=\nu\in\Delta_a$. Then
the probability that the sojourn time will exceed $t$ is
$\exp\left(\int_0^t\phi(s,\nu)\Lambda
1_{h^{-1}(a)}ds\right)$. If the following jump occurs at time $T_n<\infty$
then $Y$ jumps to $b$ with probability $q\left(\Pi^\mu_{T_{n}-},b\right)$,
which depends on the position $\Pi^\mu_{T_{n}-}= \phi (T_n-T_{n-1}, \nu)$
of the process $\Pi^\mu$ immediately before
the jump. Thus, conditionally on the ``past" $(Y_t)_{0\le t\le T_{n-1}}$,
the distribution of the ``future" $(Y_t)_{ t> T_{n-1}}$ is determined
by the present value of the filtering process $\Pi^\mu_{T_{n-1}}$.
\end{em}
\end{remark}

As a preparation for the proof of theorem \ref{mainyepi} we need the following lemma.

\begin{lemma}
Let $J_t^b$ denote the number of jumps of $Y$ to $b\in O$ in
the time interval $(0,t]$:
$$
J_t^b=\sum_{k=1}^\infty 1_{T_k\le t}1_{Y_{T_k}=b},
\qquad t\ge 0.
$$
Then the process
$$
M_t:= J_t^b-\int_0^t \Pi^\mu_s\Lambda 1_{h^{-1}(b)}\, 1_{Y_s\ne b}\,ds,
\qquad t\ge 0,
$$
is a martingale with respect to $(\caly^0_t)$.

Moreover
\begin{equation}\label{attesacondsigmannulla}
    E_\mu \left[ M_{t\wedge T_n}-M_{t\wedge T_{n-1}}\,\big| \, \Sigma_{n-1}\right]=0,
\qquad t\ge 0,\,n\ge 1.
\end{equation}
\end{lemma}

 \noindent {\bf Proof.}
 We start recalling the well known fact that for every $f:I\to\R$ the process
 $$
 \Pi^\mu_tf-\Pi^\mu_0f-\int_0^t \Pi^\mu_s\Lambda f\, ds,\qquad t\ge 0,
 $$
 is a $(\caly^0_t)$-martingale (see e.g. \cite{rowi} formula VI-(8.17)).
Choosing $f=1_{b}\circ h=1_{h^{-1}(b)}$ we have
$\Pi^\mu_t f=E_\mu [1_b(Y_t)|\caly^0_t]=1_b(Y_t)$ and we deduce  that
$$
m_t:= 1_b(Y_t)-1_b(Y_0)- \int_0^t \Pi^\mu_s\Lambda 1_{h^{-1}(b)}\, ds,\qquad t\ge 0,
 $$
 is a martingale. Next we set
 $$
 M_t:= \int_{(0,t]} 1_{I\backslash \{b\}}(Y_{s-})\, dm_s
 = \int_{(0,t]} 1_{Y_{s-}\neq b}\, dm_s,\qquad t\ge 0.
 $$
 $M_t$ is defined as a pathwise Stieltjes integral and it is a martingale, since
the integrand process is  $(\caly^0_t)$-predictable and bounded. Finally we note that
$$
\begin{array}{lll}
 M_t&=&\dis
 \int_{(0,t]} 1_{Y_{s-}\neq b}\, d1_b(Y_s)-
  \int_{(0,t]} 1_{Y_{s-}\neq b}\Pi^\mu_s\Lambda 1_{h^{-1}(b)}\, ds
 \\ &=&\dis
 J^b_t-   \int_0^t 1_{Y_{s}\neq b}\Pi^\mu_s\Lambda 1_{h^{-1}(b)}\, ds.
\end{array}
 $$
 The last equality holds since $d1_b(Y_s)$ is a measure equal to $1$ at each point
 where $Y$ jumps to $b$, and equal to $-1$ at each point where $Y$ leaves $b$,
 and since  we have $Y_{s-}=Y_{s}$ almost everywhere with respect to
 the Lebesgue measure $ds$.

 For fixed $t$, the stopped process $(M_{s\wedge t} )_{s\ge 0}$ is a uniformly
 integrable $(\caly^0_t)$-martingale and by optional stopping
 $M_{t\wedge T_{n-1}}=
E_\mu \left[ M_{t\wedge T_n}\,\big| \, \caly^0_{T_{n-1}}\right]$. Since
$\Sigma_{n-1}\subset \caly^0_{T_{n-1}}$ this
 proves the last assertion of the lemma.
 \qed

 \noindent {\bf Proof of theorem \ref{mainyepi}.}
 We continue the notation of the previous lemma and we consider,
 for $t\ge0$ and $n\ge 1$,
 $$
 M_{t\wedge T_n}-M_{t\wedge T_{n-1}}=
 J^b_{t\wedge T_n}-J^b_{t\wedge T_{n-1}}-
 \int_{t\wedge T_{n-1}}^{t\wedge T_n}\Pi^\mu_s\Lambda 1_{h^{-1}(b)}\, 1_{Y_s\ne b}\,ds.
$$
Since  $J^b_{t\wedge T_n}=\sum_{k=1}^n 1_{T_k\le t}1_{Y_{T_k}=b}$ we have
$ J^b_{t\wedge T_n}-J^b_{t\wedge T_{n-1}}=1_{T_n\le t}1_{Y_{T_n}=b}$.
Next
$$
\begin{array}{lll}\dis
 \int_{t\wedge T_{n-1}}^{t\wedge T_n}\Pi^\mu_s\Lambda 1_{h^{-1}(b)}\, 1_{Y_s\ne b}\,ds
 &=&\dis
  \int_{0}^{t}1_{ T_{n-1}\le s<T_n }\Pi^\mu_s\Lambda 1_{h^{-1}(b)}\, 1_{Y_s\ne b}\,ds
  \\&=&\dis
  \int_{0}^{t}1_{ T_{n-1}\le s<T_n }\phi(s-T_{n-1},\Pi^\mu_{T_{n-1}})\Lambda 1_{h^{-1}(b)}\,
  1_{Y_{T_{n-1}}\ne b}\,ds,
  \end{array}
  $$
  since for $T_{n-1}\le s<T_n$ we have $\Pi^\mu_s=\phi(s-T_{n-1},\Pi^\mu_{T_{n-1}})$
  and $Y_s=Y_{T_{n-1}}$.  Then from (\ref{attesacondsigmannulla}) it follows that
$$
P_\mu ({T_n\le t},{Y_{T_n}=b}\, |\, \Sigma_{n-1})=
E_\mu\left[
  \int_{0}^{t}1_{ T_{n-1}\le s<T_n }\phi(s-T_{n-1},\Pi^\mu_{T_{n-1}})\Lambda 1_{h^{-1}(b)}\,
  \,ds\, 1_{Y_{T_{n-1}}\ne b}\, \big|\, \Sigma_{n-1} \right].
$$
By standard arguments, for every $t\ge 0$ we can choose a version of the conditional probability
$P_\mu (T_n> t\, |\, \Sigma_{n-1})$ in such a way that the function
$t\to P_\mu (T_n> t\, |\, \Sigma_{n-1})$ is nonincreasing and right-continuous $P_\mu$-a.s.
In particular, $P_\mu (T_n> t\, |\, \Sigma_{n-1})(\omega)$ is jointly measurable in $(\omega,t)$.
An application of the Fubini
theorem shows that
\begin{equation}\label{congiuntacondizuno}
    P_\mu ({T_n\le t},{Y_{T_n}=b}\, |\, \Sigma_{n-1})=
  \int_{0}^{t}\!P_\mu (T_n>s|\Sigma_{n-1})\,1_{ T_{n-1}\le s }\phi(s-T_{n-1},\Pi^\mu_{T_{n-1}})
  \Lambda 1_{h^{-1}(b)}\,
   1_{Y_{T_{n-1}}\ne b}\,ds.
\end{equation}
Now we sum over all $b\in O$. Denote $a=Y_{T_{n-1}}$ and $\nu = \phi(s-T_{n-1},\Pi^\mu_{T_{n-1}})$ for short.
Then $ \Pi^\mu_{T_{n-1}}\in\Delta_a$ and consequently $\nu\in\Delta_a$ since the flow
leaves $\Delta_a$ invariant. The identity
$$\nu \Lambda 1_{h^{-1}(a)}+
\sum_{b\in O,b\neq a}\nu \Lambda 1_{h^{-1}(b)}=0,
$$
already noticed before, shows that
$$\sum_{b\in O}
\phi(s-T_{n-1},\Pi^\mu_{T_{n-1}})\Lambda 1_{h^{-1}(b)}\,
   1_{Y_{T_{n-1}}\ne b}= -\phi(s-T_{n-1},\Pi^\mu_{T_{n-1}})\Lambda 1_{h^{-1}(Y_{T_{n-1}})}
   $$
and we arrive at
\begin{equation}\label{condizdue}
    P_\mu ({T_n\le t}\, |\, \Sigma_{n-1})=
-  \int_{0}^{t}P_\mu (T_n>s|\Sigma_{n-1})\,1_{ T_{n-1}\le s }\,
\phi(s-T_{n-1},\Pi^\mu_{T_{n-1}})\Lambda 1_{h^{-1}(Y_{T_{n-1}})}
\,ds.
\end{equation}
It follows that, $P_\mu$-a.s.,  the function
$t\to P_\mu (T_n> t\, |\, \Sigma_{n-1})$ is absolutely continuous
with derivative
$$
 \frac{d}{dt}  P_\mu (T_n> t\, |\, \Sigma_{n-1})
 = P_\mu (T_n>t|\Sigma_{n-1})\,1_{ T_{n-1}\le t }\,
\phi(t-T_{n-1},\Pi^\mu_{T_{n-1}})\Lambda 1_{h^{-1}(Y_{T_{n-1}})}.
$$
Together with the condition $P_\mu (T_n\ge 0|\Sigma_{n-1})=1$ this implies
\begin{equation}\label{condiztre}
    P_\mu ({T_n> t}\, |\, \Sigma_{n-1})=
    \exp\left(
    \int_{0}^{t}1_{ T_{n-1}\le s }\,
\phi(s-T_{n-1},\Pi^\mu_{T_{n-1}})\Lambda 1_{h^{-1}(Y_{T_{n-1}})}
\,ds\right).
\end{equation}
This equality will be used in three ways. First, substituting in
(\ref{congiuntacondizuno}), we obtain a formula for the joint distribution
of $T_n$ and $Y_{T_n}$ conditional on $\Sigma_{n-1}$:
\begin{equation}\label{congiuntacondizdue}
\begin{array}{lll}\dis
    P_\mu ({T_n\le t},{Y_{T_n}=b}\, |\, \Sigma_{n-1})
    &=&\dis
  \int_{0}^{t}
  \exp\left(
    \int_{0}^{s}1_{ T_{n-1}\le r }\,
\phi(r-T_{n-1},\Pi^\mu_{T_{n-1}})\Lambda 1_{h^{-1}(Y_{T_{n-1}})}
\,dr\right)
  \cdot
  \\&&\dis \cdot
  1_{ T_{n-1}\le s }\phi(s-T_{n-1},\Pi^\mu_{T_{n-1}})
  \Lambda 1_{h^{-1}(b)}\,
   1_{Y_{T_{n-1}}\ne b}\,ds.
   \end{array}
\end{equation}
Second, (\ref{condiztre}) shows that the random variable $T_n$ (which
may take the value $\infty$) has the property that its conditional
distribution with respect to $\Sigma_{n-1}$, restricted to $[0,\infty)$,
 posses a density $d_n$ with respect to the Lebesgue measure and $d_n$
 can be computed by differentiating the right-hand
 side of (\ref{condiztre}):
\begin{equation}\label{conddeenstn}
\begin{array}{lll}
    d_n(t)&=&\dis -\exp\left(
    \int_{0}^{t}1_{ T_{n-1}\le s }\,
\phi(s-T_{n-1},\Pi^\mu_{T_{n-1}})\Lambda 1_{h^{-1}(Y_{T_{n-1}})}
\,ds\right)\cdot
  \\&&\dis \cdot
1_{ T_{n-1}\le t }\,
\phi(t-T_{n-1},\Pi^\mu_{T_{n-1}})\Lambda 1_{h^{-1}(Y_{T_{n-1}})}.
\end{array}
\end{equation}
Third, (\ref{sncondizy}) can be easily deduced from (\ref{condiztre}) as follows:
since
$$
P_\mu(S_n>t,T_{n-1}<\infty\,|\,\Sigma_{n-1})
=
P_\mu(T_n>t+T_{n-1},T_{n-1}<\infty\,|\,\Sigma_{n-1}),
$$
and
since $T_{n-1}$ is $\Sigma_{n-1}$-measurable we deduce from
(\ref{condiztre}) that
$$
\begin{array}{l}
P_\mu(S_n>t,T_{n-1}<\infty\,|\,\Sigma_{n-1})
\\
\qquad=\dis
1_{T_{n-1}<\infty}
\exp\left(\int_0^{t+T_{n-1}}1_{ T_{n-1}\le s }\,\phi(s-T_{n-1},\Pi^\mu_{T_{n-1}})\Lambda
1_{h^{-1}(Y_{T_{n-1}})}ds\right)
\\\qquad
=\dis
1_{T_{n-1}<\infty}
\exp\left(\int_0^t\phi(s,\Pi^\mu_{T_{n-1}})\Lambda
1_{h^{-1}(Y_{T_{n-1}})}ds\right).
\end{array}
$$
It remains  to prove (\ref{ytncondizy}). To this end we compute
$$
E_\mu\left[
q\left(\Pi^\mu_{T_{n}-},b\right)1_{T_n\le t}\,|\,\Sigma_{n-1} \right]
=E_\mu\left[
q\left(\phi (T_n-T_{n-1}, \Pi^\mu_{T_{n-1}}),b\right)1_{T_n\le t}\,|\,\Sigma_{n-1} \right].
$$
Noting that $T_{n-1}$ and $\Pi^\mu_{T_{n-1}}$ are $\Sigma_{n-1}$-measurable
and recalling the expression of the conditional density (\ref{conddeenstn}) we obtain
$$\begin{array}{l}\dis
E_\mu\left[
q\left(\Pi^\mu_{T_{n}-},b\right)1_{T_n\le t}\,|\,\Sigma_{n-1} \right]
\\\dis\qquad
=
\int_0^t
q\left(\phi (s-T_{n-1}, \Pi^\mu_{T_{n-1}}),b\right)\,d_n(s)\, ds
\\\dis\qquad
=-\int_0^t
q\left(\phi (s-T_{n-1}, \Pi^\mu_{T_{n-1}}),b\right)\,
\phi(s-T_{n-1},\Pi^\mu_{T_{n-1}})\Lambda 1_{h^{-1}(Y_{T_{n-1}})}
\cdot
\\ \qquad\quad\cdot \dis 1_{ T_{n-1}\le s }\,
\exp\left(
    \int_{0}^{s}1_{ T_{n-1}\le r }\,
\phi(r-T_{n-1},\Pi^\mu_{T_{n-1}})\Lambda 1_{h^{-1}(Y_{T_{n-1}})}
\,dr\right)
\, ds.
\end{array}
$$
For a moment
denote $Y_{T_{n-1}}$ by $a$ and
$\phi(s-T_{n-1},\Pi^\mu_{T_{n-1}})$ by $\nu$.
Then
$\Pi^\mu_{T_{n-1}}\in\Delta_a$, which implies
$\nu\in\Delta_a$
by the invariance
of $\Delta_a$ under the flow $\phi$.
If $a=b$ then
$q\left(\nu,b\right)=0$. If
$a\neq b$ then
$q\left(\nu,b\right)
\nu\Lambda 1_{h^{-1}(a)}
= -\nu\Lambda 1_{h^{-1}(b)}$
as noticed in (\ref{cancq}).
So we obtain
$$\begin{array}{l}\dis
E_\mu\left[
q\left(\Pi^\mu_{T_{n}-},b\right)1_{T_n\le t}\,|\,\Sigma_{n-1} \right]
\\\dis\qquad
=1_{Y_{T_{n-1}}\neq b}\int_0^t
\phi(s-T_{n-1},\Pi^\mu_{T_{n-1}})\Lambda 1_{h^{-1}(b)}
\cdot
\\ \qquad\quad
\cdot \dis 1_{ T_{n-1}\le s }\,
\exp\left(
    \int_{0}^{s}1_{ T_{n-1}\le r }\,
\phi(r-T_{n-1},\Pi^\mu_{T_{n-1}})\Lambda 1_{h^{-1}(Y_{T_{n-1}})}
\,dr\right)
\, ds.
\end{array}
$$
Comparing with (\ref{congiuntacondizdue}) we conclude that
$$
    P_\mu ({T_n\le t},{Y_{T_n}=b}\, |\, \Sigma_{n-1})=
E_\mu\left[
q\left(\Pi^\mu_{T_{n}-},b\right)1_{T_n\le t}\,|\,\Sigma_{n-1} \right].
$$
Since  $\Sigma_{n-1}^+$ is generated by $\Sigma_{n-1}$ and by $T_n$, this
immediately implies
(\ref{ytncondizy}).
\qed

\subsection{An application to exit time distributions}
\label{tempouscita}
As an example of the kind of results that can follow from
theorem \ref{mainyepi} we prove an explicit formula for
the law of the exit time of a finite Markov chain from a given set.

We start from
formula  (\ref{sncondizy}) written for $n=1$:
$$
P_\mu(S_1>t\,|\,\sigma(Y_0))=
\exp\left(\int_0^t\phi(s,\Pi^\mu_{0})\Lambda
1_{h^{-1}(Y_0)}ds\right),\qquad t\ge 0.
$$
Assume in addition that $\mu=\delta_i$ is the measure concentrated
at some $i\in I$ and let $a=h(i)\in O$.
We denote $
P_\mu$ by $P_i$ and we note that $Y_0=a$ and
$\Pi^{\mu}_{0}=\delta_i$, so we obtain
$$
P_i(S_1>t)=
\exp\left(\int_0^t\phi(s,\delta_i)\Lambda
1_{h^{-1}(a)}ds\right),\qquad t\ge 0.
$$
Finally assume that $O=\{a,b\}$ consists of exactly two points, and
denote $A=h^{-1}(a)$. Then $S_1$ coincides
with the first exit time from $A$:
\begin{equation}\label{defexitime}
    \tau=\inf \{t\ge0\,:\,X_t\notin A\}.
\end{equation}
Note that $y(t):=\phi(t,\delta_i)$ is a solution of the differential equation
$$
y'(t)= 1_{h^{-1}(a)}* (y(t)\Lambda) - (y(t)\Lambda
1_{h^{-1}(a)})\,y(t),\qquad t\ge 0,
$$
with initial condition $y(0)=\delta_i$.
By {proposition} \ref{lemmaflow} there exists a unique global solution
 with values in  $\Delta_a$. So the only possibly nonzero components
 are $y(t,j)$, for $j\in A$, and  the equation can be  written in scalar form as
\begin{equation}\label{eqscalare}
y'(t,j)= \sum_{k\in A} y(t,k)\lambda_{kj} - \left(\sum_{k,h\in A}
y(t,k)\lambda_{kh}
\right)\,y(t,j),\qquad j\in A,t\ge 0,
\end{equation}
with initial conditions
\begin{equation}\label{eqscalarecondiniz}
y(t,i)=1,\qquad
y(t,j)= 0, \; j\in A,\,j\neq i.
\end{equation}
We have finally proved the following:

\begin{proposition}\label{exitime}
Suppose $X$ is a time-homogeneous Markov chain in a finite set $I$ with
rate transition matrix $\Lambda=(\lambda_{nm})_{n,m\in I}$. Let $A$ be a proper subset of $I$,
let $P_i$ denote the law of the chain starting at $i\in A$ and let
$\tau$ be the first exit time from $A$ as defined in (\ref{defexitime}).
Then we have
$$
P_i(\tau>t)=
\exp\left(\int_0^t\sum_{k,h\in A}
y(s,k)\lambda_{kh}\;ds\right), \qquad t\ge 0,
$$
where $y(t,j)$, ($j\in A,t\ge0$) is the unique solution of
(\ref{eqscalare})- (\ref{eqscalarecondiniz}).
\end{proposition}

\begin{remark}{\em
{Proposition} \ref{exitime} is a statement on
the law of the exit time of a finite Markov chain from a given set.
Since it
is not directly related to filtering theory,
but it  rather concerns a basic topic in the theory of Markov chains,
it may be possibly proved by different arguments. However we were not
able to find a reference providing such an explicit formula.
}
\end{remark}

\begin{remark}{\em
{Proposition} \ref{exitime} follows from an application of
formula  (\ref{sncondizy}). Similar arguments based on
(\ref{ytncondizy}) provide the distribution of the Markov chain
at the exit time from the set $A$. We omit the details.
}
\end{remark}

\section{Filtering processes and piecewise-deterministic Markov processes}
\label{sectionpdp}

The main purpose of this section is to show
that the filtering process is a piecewise-deterministic Markov process
(PDP) in the sense of
Davis \cite{Da2}, \cite{Da}, and
to present some consequences. To this end
we first  recall the definition of this class of processes.

\subsection{Piecewise-deterministic Markov processes (PDPs)}

We limit ourselves
to the special case when the state space of the PDP is the effective simplex
$\Delta_e$, since this is the only case we will deal with.
Other differences from the general framework considered in
\cite{Da} are pointed out in remark \ref{specialpdprem} below.
We recall that $\Delta_e$ is a disjoint union $\cup_{a\in O}\Delta_a$ where each
$\Delta_a$ is a  compact subset of the euclidean space. We assume that we are also given
the following objects.

\begin{enumerate}
\item A {\em flow} $\phi$ on $\Delta_e$. By this we mean a continuous function
$\phi:\R_+\times \Delta_e\to \Delta_e$ such that
$\phi(t,\phi(s,x))=\phi(t+s,x)$ for $t,s\ge0$ and $x\in \Delta_e$ and
leaving each set $\Delta_a$ invariant, i.e.
$\phi(t,x)\in \Delta_a$ if $x\in \Delta_a$.

\item A {\em jump rate function} $\lambda:\Delta_e\to\R_+$. We require
that it is measurable and
that for every $x\in \Delta_e$ there exists $\epsilon>0$ (depending on $x$) such that
$\int_0^\epsilon \lambda(\phi(s,x))ds<\infty$.

\item A {\em transition measure} $Q$ on $(\Delta_e, \calb(\Delta_e))$,
i.e.  a stochastic kernel $Q(x,A)$ defined for
$x\in\Delta_e$, $A\in\calb(\Delta_e)$. We require that $Q(x,\{x\})=0$
for every $x\in \Delta_e$.

\end{enumerate}

A process $(Z_t)_{t\ge 0}$, defined on some probability space $(\Omega',\calf',P')$,
is called a PDP with respect to $(\phi,\lambda,Q)$, starting at $z_0\in\Delta_e$, if
there exists a sequence of nondecreasing random variables $T_n:\Omega'  \to [0,\infty]$
($n\ge 1$),
such that the following holds.

\begin{enumerate}
\item[(i)]
$Z_0=z_0$ and
the trajectories of $(Z_t)$ are cadlag functions, with discontinuities occurring
precisely at times $(T_n)$, $P'$-a.s.

\item[(ii)] Among jump times the process evolves deterministically along the flow; more
precisely we have,  for $n\ge1$, $P'$-a.s.,
 \begin{equation}    \label{trajpdp}
Z_t=\phi(t-T_{n-1}, Z_{T_{n-1}}) {\rm    \quad for\quad}  T_{n-1}\le t <T_n,
\end{equation}
where we set $T_0=0$.
This implies that the random variables
$Z_{T_n-}:=\lim_{t\to {T_n}, t<{T_n}}Z_t$,
defined on $\{T_n<\infty\}$,  are also given by the formula
$$
Z_{T_n-}= \phi(T_n-T_{n-1},  Z_{T_{n-1}})= \phi(S_n,  Z_{T_{n-1}}),
\qquad {\rm on\;} \{T_n<\infty\},\; n\ge 1,
$$
where $S_n=T_n-T_{n-1}$ ($S_n$ are defined for $n\ge 1$ on the event $\{T_{n-1}<\infty\}$).

\item[(iii)]
For $t>0$ and $A\in \calb(\Delta_e)$,
\begin{equation}    \label{inizpdp}
\begin{array}{l}\dis P'(T_1>t)
=
\exp\left(-\int_0^t\lambda(\phi(s,z_0))\,ds\right),
\\\dis
P'(Z_{T_1}\in A\,|\,
T_1)
=Q(Z_{T_1-}, A)= Q(\phi(T_1,z_0), A),\qquad {\rm on\;} \{T_1<\infty\},
\end{array}
\end{equation}
and for $n\ge2$,
\begin{equation}    \label{jointpdp}
\begin{array}{l}\dis
P'(S_n>t\,|\,
T_1,Z_{T_1},\ldots, T_{n-1},Z_{T_{n-1}})
=
\exp\left(-\int_0^t\lambda(\phi(s,Z_{T_{n-1}}))
\,ds\right)
\; {\rm on\;} \{T_{n-1}<\infty\}
\\\dis P'(Z_{T_n}\in A\,|\,
T_1,Z_{T_1},\ldots, T_{n-1},Z_{T_{n-1}},T_n )
=Q(Z_{T_n-}, A)
\\\dis\qquad\qquad\qquad\qquad\qquad\qquad\qquad
\qquad\qquad\,
= Q(\phi(S_n,  Z_{T_{n-1}}), A)
\; {\rm on\;} \{T_n<\infty\}.
\end{array}
\end{equation}

\end{enumerate}

Formulae (\ref{inizpdp})-(\ref{jointpdp}) allow to interpret
$Q(x,A)$ as the probability to find the process in the set $A\subset \Delta_e$
at a jump time, conditional to the fact that the process was in $x\in\Delta_e$
 immediately before the jump. They also explain the terminology
 {\em jump rate} for the function $\lambda$.

Formulae (\ref{inizpdp})-(\ref{jointpdp}) also show that
$(\phi,\lambda,Q)$ and the starting point $z_0$
uniquely determine the
 finite-dimensional
distributions
of the stochastic process
$\{T_1,$ $ Z_{T_1},$ $ T_2, $ $Z_{T_1},$ $\ldots\}$.
In view of (\ref{trajpdp}) we conclude that the law of $(Z_t)$ is completely
determined by $(\phi,\lambda,Q)$ and $z_0$.

\begin{remark}\begin{em}\label{specialpdprem}
The present definition differs from the one given in \cite{Da} for the following
reasons.

\begin{enumerate}
\item[(i)]
In \cite{Da} a specific probability space is chosen, namely the one consisting
of a countable product of unit intervals equipped with the product Lebesgue measure.
However the law of the constructed process is the same, and this is what matters
in the following.

\item[(ii)] In \cite{Da} the state space $\Delta_e=\cup_{a\in O}\Delta_a$
is replaced by a general, finite or countable, union $E=\cup_{\nu}E_\nu$ of {\em open}
sets $E_\nu$ of the euclidean space. The fact that $\Delta_a$ are compact
does not affect the basic results. Similar slight differences
are sometimes present
in the literature, for instance in \cite{Ga} the state space $E$ is assumed
to be closed; also compare the discussion in \cite{Da}, at the beginning
of section 24, on possible
generalizations to cases where each $E_\nu$ may be a differentiable manifold.

\item[(iii)] In \cite{Da}, instead of the flow, a vector field is chosen as
a starting point. This field is assumed to be locally Lipschitz and to generate
a flow which is defined up to the time when it hits the boundary of $E_\nu$,
 for every starting point in $E_\nu$. In the specific situation of the filtering process
 which we are about to study we will also exhibit the vector field associated
 to the flow.

\item[(iv)] The main difference  is the fact that in \cite{Da} the trajectories of
the PDP process
$(Z_t)$ are required to jump at each time when they hit the boundary
of some $E_\nu$. Jumps at the boundary will not occur for the filtering process
presented later.
This difference does not affect the basic results we are going to use,
and in fact it results in a simplification: for instance the delicate ``boundary
conditions'' presented in \cite{Da} in connection with a description of
the extended infinitesimal generator of the PDP process are not needed.
Again, similar differences
are already  present
in the literature, for instance in \cite{Ga} a jump occurs whenever
the process hits some prescribed closed set $\Gamma\subset E$, not
necessarily the boundary of $E$.

\end{enumerate}

\end{em}
  \end{remark}

\subsection{The filtering process as a PDP}

Now we come back to the canonical set-up introduced in subsection
\ref{canonical}, we fix $\nu\in\Delta_e$ and we construct
the probability  space $(\Omega,\calf^0,P_\nu)$. In the rest of this
section all stochastic processes
will be considered under $P_\nu$, so that in particular
 $\nu$ is the initial distribution of $(X_t)$. Since $\nu$ belongs to
 the effective simplex,
as noted in remark \ref{piepinu}, the filtering process  $(\Pi^\nu_t)$
starts at $\nu$, i.e.
$\Pi^\nu_0=\nu$.
 It is our purpose to show
that $(\Pi^\nu_t)$ is a PDP and to describe explicitly the corresponding
triple $(\phi,\lambda,Q)$.

\begin{enumerate}
\item As the flow $\phi$ we take    the flow introduced in subsection \ref{simplexflow}
after proposition \ref{lemmaflow}.
We recall that it is the  flow associated to the vector field $F$ defined on $\Delta_e=
\cup_{a\in O}\Delta_a$ by the formula
$$
F(\nu)=1_{h^{-1}(a)}\cdot (\nu\Lambda)
- (\nu\Lambda 1_{h^{-1}(a)})\nu,\qquad \nu\in \Delta_a.
$$

\item As the jump rate function we take the function $\lambda:\Delta_e\to\R_+$
defined by
\begin{equation}\label{defjumprate}
\lambda(\nu )=
-\nu \Lambda
1_{h^{-1}(a)},\qquad \nu \in\Delta_a.
\end{equation}
It was shown in section \ref{sectionobservation} that
 $\lambda(\nu )\ge0$.

\item The transition measure   $Q(\nu ,A)$
 is defined for $\nu \in\Delta_e$ and $ A\in\calb(\Delta_e)$ by
\begin{equation}\label{deftransmeas}
Q(\nu ,A)=
\sum_{b\in O}
1_A( H_b[\nu \Lambda]) \,q(\nu,b),
\end{equation}
where $q(\nu,b)$ was introduced in section \ref{sectionobservation}
(we still exclude the trivial case where $h$ is constant).
Thus,
for
$\nu \in\Delta_a$,
 $Q(\nu ,\cdot)$
 is  the measure concentrated on the finite set
 $$\{H_b[\nu \Lambda]\,:\, b\in O\backslash\{ a\}\}\subset
 \Delta_e\backslash \Delta_a,
 $$ and each  point $H_b[\nu \Lambda]$ has mass
 $q(\nu,b)$. More explicitly, for
$\nu \in\Delta_a$,
$$Q(\nu ,\{H_b[\nu \Lambda]\})=
\frac{\nu \Lambda 1_{h^{-1}(b)}
}{-\nu \Lambda 1_{h^{-1}(a)}
}\, 1_{b\neq a},
$$
provided $\nu \Lambda 1_{h^{-1}(a)}\neq 0$.

\end{enumerate}

\begin{remark}\begin{em}\label{qnonfeller}
In the literature on PDPs the requirement that $Q$ is a Feller kernel
is often formulated among the assumptions:
this means that for every bounded continuous
$g:\Delta_e\to\R$, the function $\nu\mapsto \int_{\Delta_e}g(\rho)
 Q(\nu,d\rho)$ is continuous  (and necessarily bounded).
This  condition fails in general in our case, due to the occurrency
of the exceptional set where $\nu \Lambda 1_{h^{-1}(a)}= 0$ ($\nu\in\Delta_a$)
in the formulae above. However, in our case
we have the following weaker form of the Feller property of $Q$.

\end{em}
  \end{remark}

\begin{proposition}\label{quasifeller}
The function
$\nu\mapsto  \lambda(\nu)\int_{\Delta_e}g(\rho)
 Q(\nu,d\rho)$ is continuous (and obviously bounded) on $\Delta_e$
for every
bounded continuous function
$g:\Delta_e\to\R$.
\end{proposition}

 \noindent {\bf Proof.} We start from formula
(\ref{cancq}): by the previous definitions it can be
written, for $\nu \in\Delta_a$ and $b\neq a$:
$$\lambda(\nu) \,q(\nu, b)=
\nu \Lambda 1_{h^{-1}(b)}.
$$
Therefore, since $q(\nu,a)=0$ for $\nu\in\Delta_a$, we obtain
$$
    \lambda (\nu ) Q(\nu ,A)=
\sum_{b\in O,b\neq a}
1_A( H_b[\nu \Lambda]) \,
\nu \Lambda 1_{h^{-1}(b)},
\qquad
\nu \in\Delta_a,\,
A\in\calb(\Delta_e),
$$
and it follows that, for
 $\nu\in\Delta_a$,
$$
 \lambda(\nu)
 \int_{\Delta_e}g(\rho)
 Q(\nu,d\rho)=
\sum_{b\in O,b\neq a}
g( H_b[\nu \Lambda]) \,
\nu \Lambda 1_{h^{-1}(b)}.
 $$
 Now if $\nu_n\to\nu$ then for all large $n$ we have  $\nu_n\in\Delta_a$ and
we distinguish two cases.

 If $\nu \Lambda 1_{h^{-1}(a)}\neq0$ then  for $b\neq a$
 $$
H_b[\nu_n\Lambda]=
\frac{1}{\nu_n \Lambda 1_{h^{-1}(a)}}1_{h^{-1}(a)}\ast (\nu_n \Lambda)
\to
\frac{1}{\nu \Lambda 1_{h^{-1}(a)}}1_{h^{-1}(a)}\ast (\nu \Lambda)
=H_b[\nu \Lambda],
$$
and the result follows from the continuity of $g$;

If $\nu \Lambda 1_{h^{-1}(a)}=-\lambda(\nu)=0$ then for all $b\neq a$ we have
$\nu_n \Lambda 1_{h^{-1}(b)}\to \nu \Lambda 1_{h^{-1}(b)}=0$,
 which implies
$g( H_b[\nu_n \Lambda]) \,
\nu_n \Lambda 1_{h^{-1}(b)}\to 0
$ by the boundedness of $g$.
\qed

We are now ready to present the main result of this section.

\begin{theorem}\label{mainpdp}
For every $\nu\in\Delta_e$  the filtering process  $(\Pi^\nu_t)$,
defined in
the probability  space $(\Omega,\calf^0,P_\nu)$ and taking values
in $\Delta_e$,
is a piecewise-deterministic Markov process with respect
to the
triple $(\phi,\lambda,Q)$ defined above and with starting point $\nu$.
\end{theorem}

 \noindent {\bf Proof.}
Let $(T_n)_{n\ge 1}$ be the  sequence of  jumps times of $Y$,
with the convention that $T_n=\infty$ for all $n$ if no jump occurs
and $T_{n}<T_{n+1}=T_{n+2}=\ldots=\infty$ if precisely $n$ jumps occur. We let
 $T_0=0$ and define $S_n=T_n-T_{n-1}$
  for $n\ge 1$ on the event $\{T_{n-1}<\infty\}$.

According to equation (\ref{eqdefpifiltro}),
  $(\Pi^\nu_t)$ can be explicitly described as follows: the starting point is
$\Pi^\nu_0=\nu$ and
 for $n\ge1$,
 $$\begin{array}{ll}
\Pi^\nu_t=\phi(t-T_{n-1}, \Pi^\nu_{T_{n-1}})& {\rm    \;for\;}  T_{n-1}\le t <T_n,
\\
\Pi^\nu_{T_n-}=\phi(S_n, \Pi^\nu_{T_{n-1}}), &
\Pi^\nu_{T_n}=H_{Y_{T_n}}[\Pi^\nu_{T_n-}\Lambda].
\end{array}
$$
It is clear that $\Pi^\nu$ can jump only at times when $Y$ jumps.
We now claim that each $T_n$ ($n\geq 1$)
is also a jump time of $\Pi^\nu$ and therefore the jump times of $\Pi^\nu$
and $Y$ coincide. Indeed,  if $a$ denotes $Y_{T_n-}=Y_{T_{n-1}}$ and $b$ denotes
$Y_{T_n}$ then
$a\neq b$ and since
$\Pi^\nu_{T_n-}\in \Delta_a$, $\Pi^\nu_{T_n}\in \Delta_b$,
it follows  that $\Pi^\nu_{T_n-}
\neq \Pi^\nu_{T_n}$.

In (\ref{defsigma}) we defined for $n\ge 1$
$$
    \Sigma_{n-1}=\sigma (Y_0,T_1,Y_{T_1},\ldots, T_{n-1},Y_{T_{n-1}}),\quad
\Sigma^+_{n-1}=\sigma (Y_0,T_1,Y_{T_1},\ldots, T_{n-1},Y_{T_{n-1}},T_n ).
$$
Note that $Y_0$ is constant $P_\nu$-a.s. since if $\nu\in\Delta_a$ for some $a\in O$ then
$Y_0=a$. Also note that  by the filtering equation
(\ref{eqdefpifiltro})
$T_1,Y_{T_1},\ldots, T_{n-1},Y_{T_{n-1}}$
uniquely determine $T_1,\Pi^\nu_{T_1},\ldots, T_{n-1},\Pi^\nu_{T_{n-1}}$. However,
the converse is also true, since if $\Pi^\nu_{T_k}\in\Delta_a$
for some $a\in O$ at a jump time $T_k<\infty$ then
$Y_{T_k}=a$. So we conclude that up to $P_\nu$-null sets,  for $n\ge 1$,
$$
    \Sigma_{n-1}=\sigma (T_1,\Pi^\nu_{T_1},\ldots, T_{n-1},\Pi^\nu_{T_{n-1}}),\quad
\Sigma^+_{n-1}=\sigma (T_1,\Pi^\nu_{T_1},\ldots, T_{n-1},\Pi^\nu_{T_{n-1}},T_n ),
$$
and in particular $\Sigma_0$ is the trivial $\sigma$-algebra and $\Sigma_0^+=\sigma(T_1)$.

Formula (\ref{sncondizy}) in
theorem \ref{mainyepi}
shows that
for $t\ge 0$ and  $n\ge1$,
$$
P_\nu(S_n>t,T_{n-1}<\infty\,|\,\Sigma_{n-1})=
\exp\left(\int_0^t\phi(s,\Pi^\nu_{T_{n-1}})\Lambda
1_{h^{-1}(Y_{T_{n-1}})}ds\right)\,1_{T_{n-1}<\infty}.
$$
Denoting $Y_{T_{n-1}}$ by $a$, we have
$\Pi^\nu_{T_{n-1}}\in\Delta_a$ and
$\phi(s,\Pi^\nu_{T_{n-1}})\in\Delta_a$ by the invariance of $\Delta_a$ under
the flow. So by the definition of the jump rate function $\lambda$
(formula (\ref{defjumprate}))
we have
$\phi(s,\Pi^\nu_{T_{n-1}})\Lambda
1_{h^{-1}(Y_{T_{n-1}})}=-\lambda (\phi(s,\Pi^\nu_{T_{n-1}}))$ and we conclude
that on the event $\{T_{n-1}<\infty\}$ we have
\begin{equation}\label{sojournpdp}
    P_\nu(S_n>t\,|\,\Sigma_{n-1})=
\exp\left(-\int_0^t\lambda (\phi(s,\Pi^\nu_{T_{n-1}}))\, ds\right).
\end{equation}
This formula shows that the sojourn times $S_n$ have the required
conditional distributions.

Now we proceed to compute the conditional distributions of $\Pi^\nu_{T_n}$.
We
 note that at each jump time $T_n< \infty$ ($n\geq 1$) we have
$ \Pi^\nu_{T_n}=H_{Y_{T_n}}[\Pi^\nu_{T_n-}\Lambda]$ and therefore
$$
\{ \Pi^\nu_{T_n}\in A,T_n< \infty\}
=\bigcup_{b\in O} \{H_b[\Pi^\nu_{T_n-}\Lambda]\in A,
 Y_{T_n}=b, T_n< \infty\}.
$$
Since
$\Pi^\nu_{T_n-}=\phi(S_n, \Pi^\nu_{T_{n-1}})$
on $\{T_n< \infty\}$ the set
$\{ H_b[\Pi^\nu_{T_n-}\Lambda]\in A,T_n< \infty\}$
belongs to $\Sigma^+_{n-1}$
 and
therefore
$$
P_\nu (\Pi^\nu_{T_n}\in A,T_n< \infty \,|\,
\Sigma^+_{n-1})
= \sum_{b\in O}
1_{T_n< \infty}
1_A( H_b[\Pi^\nu_{T_n-}\Lambda]) \,
P_\nu (Y_{T_n}=b,T_n< \infty \,|\,
\Sigma^+_{n-1}).
$$
The right-hand side can be computed using formula (\ref{ytncondizy}) in
theorem \ref{mainyepi} and we obtain
$$
P_\nu (\Pi^\nu_{T_n}\in A,T_n< \infty \,|\,
\Sigma^+_{n-1})
= \sum_{b\in O}
1_{T_n< \infty}
1_A( H_b[\Pi^\nu_{T_n-}\Lambda]) \,
q\left(\Pi^\nu_{T_{n}-},b\right)
= 1_{T_n< \infty} \, Q(\Pi^\nu_{T_n-},A),
$$
where the last equality follows from the definition of the transition measure $Q$
(formula (\ref{deftransmeas})).

Together with  (\ref{sojournpdp}), this shows that the properties
 (\ref{inizpdp}) and  (\ref{jointpdp})
hold true and therefore
 $(\Pi^\nu_t)$ is a PDP with the prescribed
jump rate function
$\lambda$ and transition measure
$Q$.
\qed

At this point several properties of the filtering process might be
stated as immediate consequences of general results on PDPs. For
instance, an explicit description of the extended generator
of its transition semigroup can be given in term of the
triple $(\phi,\lambda,Q)$: see \cite{Da}. In
section \ref{optimalstopping} below we will use known results
on standard optimal stopping problems for PDPs in order to solve
an optimal stopping problem with partial observation for the process $X$.
In the present section we will exploit general knowledge
on PDPs in order to prove that the filtering process
has the Feller property.

\subsection{The Feller property of the filtering process}

We still use the canonical set-up of subsection
\ref{canonical} and
for every $\nu\in\Delta_e$ we consider again the filtering process  $(\Pi^\nu_t)$
defined in
the probability  space $(\Omega,\calf^0,P_\nu)$ and taking values
in $\Delta_e$. It
is a piecewise-deterministic Markov process with respect
to the
triple $(\phi,\lambda,Q)$ defined above and with starting point $\nu$.
Here we wish to investigate further properties of its transition semigroup
$(R_t)$ introduced before
proposition \ref{filtromarkoviano}.
We denote by $C(\Delta_e)$ (respectively, $B(\Delta_e)$) the space of real
continuous (respectively, Borel measurable) functions on $\Delta_e$.
Since $\Delta_e$ is compact, we have $C(\Delta_e)\subset B(\Delta_e)$.
It is clear that $R_t$ maps $B(\Delta_e)$ into $B(\Delta_e)$ for all $t\ge0$.
The fact that $R_t$ also maps $C(\Delta_e)$ into $C(\Delta_e)$ is known as
Feller property and it is proved in the following proposition, together with
a statement on continuity which shows
 that  $(R_t)$ is a strongly continuous semigroup on the space
$C(\Delta_e)$ equipped with the supremum norm.

\begin{proposition}\label{pdpfeller}
For every $f\in C(\Delta_e)$ we have
$$
R_tf\in C(\Delta_e), \qquad t\ge 0,
$$
and
$
R_tf\to f
$ uniformly on $\Delta_e$ as $t\to 0$.
\end{proposition}

 \noindent {\bf Proof.}
 Let us define an operator $G$
 acting on continuous and bounded functions
 $\psi:\R_+\times \Delta_e\to\R$ by the formula
 $$
 G\psi(t,\nu)=f(\phi(t,\nu))e^{-M(t,\nu)}
 +
 \int_0^t\int_{\Delta_e}\psi(t-s,\rho)\,Q(\phi(s,\nu),d\rho)\, \lambda(\phi(s,\nu))
 e^{-M(s,\nu)}ds,
 $$
 for $ t\ge 0$, $\nu\in\Delta_e$,
 where
 $M(t,\nu)=\int_0^t\lambda(\phi(s,\nu))\,ds$. Let
  $G^n$ denote the $n$-th iterate of  $G$.
 Since the function $\nu\mapsto \int_{\Delta_e}g(\rho)
 Q(\nu,d\rho)\, \lambda(\nu)$ is continuous for every $g\in C(\Delta_e)$
by proposition \ref{quasifeller}, therefore
 $G\psi$ (and hence $G^n\psi$) is continuous and bounded on $\R_+\times \Delta_e$.
It is proved in \cite{Da}, proof of Theorem 27.6, that  for every fixed $\psi$ we have
$ G^n\psi(t,\nu)\to R_tf(\nu)$ uniformly in $\nu\in\Delta_e$ for all $t\ge 0$ as $n\to\infty$,
 so
we immediately conclude that $R_tf\in C(\Delta_e)$ for every $t\ge 0$.

To prove the final statement of the proposition note that
$$
\begin{array}{lll}
R_tf(\nu)&=& E_\nu[f(\Pi^\nu_t)1_{t<T_1}]+E_\nu[f(\Pi^\nu_t)1_{t\ge T_1}]
\\&=&
f(\phi(t,\nu))P_\nu(t<T_1)+E_\nu[f(\Pi^\nu_t)1_{t\ge T_1}],
\end{array}
$$
so that
$$
\begin{array}{lll}
|R_tf(\nu)-f(\nu)|&=& |
f(\phi(t,\nu))-f(\nu)-f(\phi(t,\nu))P_\nu(t\ge T_1)+E_\nu[f(\Pi^\nu_t)1_{t\ge T_1}]|
\\
&\le &
|f(\phi(t,\nu))-f(\nu)| +2 (\sup |f|) \,P_\nu(t\ge T_1).
\end{array}
$$
By the properties of the flow we have $|f(\phi(t,\nu))-f(\nu)|\to 0$ as $t\to 0$, uniformly
in $\nu\in\Delta_e$. Recalling formula (\ref{inizpdp}) for the distribution of
$T_1$ and denoting by $\bar\lambda$
an upper bound for the jump rate function $\lambda$ we have
$$P_\nu(T_1\ge t)=1-
\exp\left(-\int_0^t\lambda(\phi(s,\nu))\,ds\right)
\le 1-
\exp\left(- t\bar \lambda\right),
$$
which shows that
$P_\nu(T_1\ge t)\to 0$
as $t\to 0$, uniformly
in $\nu\in\Delta_e$.
  \qed

\subsection{A canonical version of the PDP filtering process}\label{PDPcanonical}

The introduction of a canonical version is useful for applications and will be used
in section \ref{optimalstopping}.
We follow the notation of  \cite{Da}, section 25. Proposition
\ref{leggefiltro}  is the only new result in this subsection, which will be applied in
section \ref{optimalstopping}.

\begin{enumerate}
\item Let $\bar\Omega$  be the set of cadlag functions $\bar\omega:\R_+\to \Delta_e$.
We denote $\bar\Pi_t(\bar\omega)=\bar\omega(t)$ for $\bar\omega\in \bar\Omega$ and $t\ge 0$,
and we introduce the $\sigma$-algebras
$$
\bar\calf^0_t=\sigma(\bar\Pi_s\,:\,s\in[0,t]),
 \qquad
\bar\calf^0=\sigma(\bar\Pi_s\,:\,s\ge0).
$$
 $(\bar\calf^0_t)_{t\ge0}$ is thus the natural filtration of $(\bar\Pi_t)_{t\ge0}$.

\item For every  $\nu\in\Delta_e$, we denote by $\bar P_\nu$ the law
of the process $(\Pi_t^\nu)$ defined in $(\Omega,\calf^0,P_\nu)$.
Thus, $\bar P_\nu$ is the probability measure on
$(\bar\Omega,\bar\calf^0)$ such that
$\bar P_\nu (\Gamma)=P_\nu (\omega\in\Omega\,:\, \Pi_\cdot^\nu(\omega)\in\Gamma)$
and $\nu$ is the starting point of $(\bar\Pi_t)$ under $\bar P_\nu$.
This definition is meaningful since the trajectory
$t\mapsto \Pi_t^\nu(\omega)$, denoted
$\Pi_\cdot^\nu(\omega)$, belongs to $\bar\Omega$ for $P_\nu$-almost
all $\omega$ and
the map $\omega\mapsto \Pi_\cdot^\nu(\omega)$
is  measurable from $(\Omega,\calf^0)$ to $(\bar\Omega,\bar\calf^0)$.

For every Borel  probability measure $Q$ on $\Delta_e$, we define a
probability $\bar P_Q$ on $(\bar\Omega,\bar\calf^0)$
by
$\bar P_Q(\Gamma)=\int_{\Delta_e}\bar P_\nu(\Gamma)\,Q(d\nu)$
for $\Gamma \in \bar\calf^0$.
Thus, $Q$ is the initial distribution of
 $(\bar\Pi_t)$ under $\bar P_Q$.

\item
We denote by $\bar\calf^Q$ the $\bar P_Q$-completion
of $\bar\calf^0$ and we assume that $\bar P_Q$ is extended to $\bar\calf^Q$ in the natural way.
We denote by $\bar\caln^Q$ the family of elements of $\bar\calf^Q$ with zero $\bar P_Q$-probability
and  we define
$$
\bar \calf^Q_t=\sigma(\bar\calf^0_t,\bar\caln^Q),\quad
\bar \calf_t=\bigcap_Q\bar \calf^Q_t,
\qquad
t\ge 0,
$$
where the intersection is taken over all Borel probability measures $Q$ on $\Delta_e$.

 $(\bar \calf_t)_{t\ge0}$ is called the natural completed filtration of $(\bar\Pi_t)$.
 It is right-continuous, i.e. for all $t\ge0$ we have
 $\bar\calf_t=\bar\calf_{t+}:=\cap_{\epsilon>0}\bar\calf_{t+\epsilon}$:
 see \cite{Da}, theorem 25.3.

\end{enumerate}

We will refer to  the PDP on $\Delta_e$
constructed above as the
$5$-tuple $(\bar\Omega,\bar\calf$, $(\bar\Pi_t)_{t\geq 0}$,
$(\bar P_\nu)_{\nu\in \Delta_e}$,
$(\bar\calf_t)_{t\geq 0})$.

Now suppose that $\mu\in\Delta$ and consider the probability
space $(\Omega,\calf^0,P_\mu)$ and  the filtering process
$(\Pi^\mu_t)$ defined by (\ref{eqdefpifiltro}).   Since
$\mu$ does not belong to $\Delta_e$ in general, the law of
$(\Pi^\mu_t)$ is not immediately known.

\begin{proposition} \label{leggefiltro} For every $\mu\in\Delta$ the law of
$(\Pi^\mu_t)$ under $P_\mu$  is $\bar P_Q$, where $Q$ is
the Borel probability measure on $\Delta_e$
concentrated at points
 $H_a[\mu]\in\Delta_e$ ($a\in O$) such that
 $$
 Q(\{H_a[\mu]\})=\mu(h^{-1}(a)),\qquad a\in O.
 $$
\end{proposition}

 \noindent {\bf Proof.}
 Define $\rho\in\Delta$ setting $\rho=\sum_{a\in O}\mu(h^{-1}(a))\,H_a[\mu]$.
  Since $\rho$ is a convex combination of probabilities
 $H_a[\mu]$ in $\Delta_e$, the law of $(\Pi^\mu_t)$ under $P_\rho$
 is $\sum_{a\in O}\mu(h^{-1}(a))\,\bar P_{H_a[\mu]}=\bar P_Q$.
 So to finish the proof it is enough to show that the laws of $(\Pi^\mu_t)$ under $P_\mu$ and $P_\rho$
 are the same. The filtering equation (\ref{eqdefpifiltro}) and the definition of the
 observation process $Y_t=h(X_t)$ imply that the trajectories  $\Pi^\mu_\cdot (\omega)$
 of the filtering process
are a deterministic functional of the corresponding trajectories  $X_\cdot (\omega)$ of
the process $X$. So it is enough to show that $P_\mu$ and $P_\rho$ coincide
on the canonical space $(\Omega,\calf^0)$. Since the generator $\Lambda$ is fixed, it only remains
to check that the law of $X_0$ is the same under $P_\mu$ and $P_\rho$.
This is, however, immediate, since recalling the definition of the operator $H$ in
subsection \ref{sectionoperatorh} we verify that
if $i\in I$ and $h(i)=b\in O$ then we have
$$P_\rho(X_0=i)= \mu(h^{-1}(b))\,H_b[\mu](i)=\mu(i)=P_\mu(X_0=i).  \qed
$$

\section{Optimal stopping with partial observation}
\label{optimalstopping}

We assume that $I$, $\Lambda$, $h$ are given as in the previous
paragraphs and we consider again the canonical set-up
introduced in subsection
\ref{canonical}. Thus $(X_t)$ is the canonical coordinate
process in the space  $\Omega$   of cadlag functions $\omega:\R_+\to I$,
$(\calf^0_t)$  is the natural filtration and $\calf^0$ is the $\sigma$-algebra
generated by $X$.
 For $\mu\in\Delta$,  $P_\mu$  denotes the probability on
$(\Omega,\calf^0)$ that makes $(X_t)$ a Markov process on $I$ with generator
$\Lambda$ and initial distribution $\mu$. The
observation process
and  its
natural filtration are still denoted $(Y_t)$ and  $(\caly^0_t)$
respectively.

For $\mu\in\Delta$ we denote by $\calf^\mu$ the $P_\mu$-completion
of $\calf^0$ and we assume that $P_\mu$ is extended to $\calf^\mu$ in the natural way.
We denote by $\caln^\mu$ the family of elements of $\calf^\mu$ with zero $P_\mu$-probability
and  we define
$$
\caly^\mu_t=\sigma(\caly^0_t,\caln^\mu),\qquad
t\ge 0.
$$
 $(\caly^\mu_t)_{t\ge0}$ is called the natural completed filtration of $(Y_t)$.
 As a consequence of the fact that $(Y_t)$ has piecewise-constant, right-continuous
 trajectories under $P_\mu$, with a finite number of jumps in every bounded interval a.s.,
it can be proved
 that
 $(\caly^\mu_t)$
  is right-continuous, i.e. for all $t\ge0$ we have
 $\caly^\mu_t=\caly^\mu_{t+}:=\cap_{\epsilon>0}\caly^\mu_{t+\epsilon}$:
 see \cite{Br}, Appendix A2, or \cite{Da}, Appendix A2. Consequently the filtered
 probability space $(\Omega,\calf^\mu,(\caly_t^\mu),P_\mu)$ satisfies
 the usual conditions.

 We denote by $\calt^\mu$ the class of functions $\tau:\Omega\to[0,\infty]$ that
 are stopping times with respect to  $(\caly^\mu_t)$.

In addition to $I$, $\Lambda$, $h$, the stopping problem with partial observation
is defined by a pair of functions
$g,l:I\to \R$, called stopping cost and running cost respectively,
 and a real number $\alpha>0$, called discount factor.
 This terminology is justified by the introduction of the following cost functional:
 $$
 J(\mu,\tau)=E_\mu\left[e^{-\alpha \tau}g(X_\tau)+\int_0^\tau
 e^{-\alpha s}l(X_s)\;ds\right],
 \qquad \mu\in\Delta,\; \tau\in\calt^\mu,
 $$
 that one tries to minimize with respect to $\tau\in\calt^\mu$.
 Here we adopt the convention that
 $e^{-\alpha \tau}g(X_\tau)=0$ if $\tau=\infty$; similar conventions
 will be tacitly used in the following.
The corresponding value function is defined by
 $$
V(\mu)=\inf_{\tau\in\calt^\mu} J(\mu,\tau),
 \qquad \mu\in\Delta.
 $$
A stopping time $\tau^{*,\mu}\in \calt^\mu$ is called optimal
(relatively to $\mu\in\Delta$) if $V(\tau^{*,\mu})=J(\mu,\tau^{*,\mu})$.
 The optimal stopping problem consists in finding characterizations
 of $V$ and giving conditions ensuring the existence of an optimal stopping time
 and a description of it, for all $\mu\in\Delta$.

\begin{remark}\label{filtrarresto}\begin{em}
In the definition of $J$ and $V$ we could replace the class $\calt^\mu$
of stopping times relatively to $(\caly^\mu_t)$ by the class
of stopping times relatively to the natural, uncompleted filtration
$(\caly^0_t)$. However, the former is much larger and, due to the
fact that $(\caly^\mu_t)$  satisfies the usual conditions,
it includes many interesting random times (for instance the first entry
time of a cadlag adapted process in a Borel set). For this
reason we have chosen the formulation above.
\end{em}
\end{remark}

For every $\mu\in\Delta$ we still denote by $(\Pi^\mu_t)$ the filtering process
defined by equation (\ref{eqdefpifiltro}).

\begin{lemma}\label{passaggiotau}
 For every $\mu\in\Delta$ and $\tau\in\calt^\mu$ we have
$$
 J(\mu,\tau)=E_\mu\left[e^{-\alpha \tau}\Pi^\mu_\tau g+\int_0^\tau
 e^{-\alpha s}\Pi^\mu_s l\;ds\right].
 $$
\end{lemma}

 \noindent {\bf Proof.} This is a direct consequence of the properties of the
 optional projections of $(e^{-\alpha t}g(X_t))$ and $(e^{-\alpha t}l(X_t))$,
 but it can also be easily proved
 by elementary considerations as follows.

 Assume first that  $\tau$ has a finite number of values.
  Excluding the trivial case of $\tau$ being constant, there exist an integer $n>1$,
  numbers $0\le t_1<\ldots<t_{n-1}<t_n=\infty$ and disjoint sets $A_i\subset \Omega$
  ($1\le i\le n$) such that $\cup_{i=1}^n A_i= \Omega$, $A_i\in \caly_{t_i}^\mu$ for
  $1\le i< n$ and $\tau= \sum_{i=1}^n t_i1_{A_i}$.
Since $\Pi^\mu_t(i)=P_\mu(X_t=i|\caly_t^\mu)$ we have
$E_\mu[g(X_{t_i})\,1_{A_i}]=E_\mu[\Pi^\mu_{t_i}g\,1_{A_i}]$ and
it follows that
$$
E_\mu[e^{-\alpha \tau}g(X_{\tau})]
=\sum_{i=1}^{n-1}E_\mu[e^{-\alpha t_i}g(X_{t_i})\,1_{A_i}]
=\sum_{i=1}^{n-1}E_\mu[e^{-\alpha t_i}\Pi^\mu_{t_i}g\,1_{A_i}]=
E_\mu[e^{-\alpha \tau}\Pi^\mu_{\tau}g].
$$

Given a general $\tau\in\calt^\mu$, let $\tau_n$ be a nonincreasing sequence of
$(\calf^\mu_{t})$-stopping times
such that $\tau_n\to\tau$
and each $\tau_n$ has a finite number of values. In the equality
$
E_\mu[e^{-\alpha \tau_n}g(X_{\tau_n})]=
E_\mu[e^{-\alpha \tau_n}\Pi^\mu_{\tau_n}g]
$
we let $n\to\infty$. We note that $X_{\tau_n}\to X_{\tau}$,
$\Pi^\mu_{\tau_n}\to\Pi^\mu_{\tau}$ by the right-continuity of
$X$ and $\Pi^\mu$, and we conclude that
$
E_\mu[e^{-\alpha \tau}g(X_{\tau})]=
E_\mu[e^{-\alpha \tau}\Pi^\mu_{\tau}g].
$

The equality
$
E_\mu\int_0^\tau
 e^{-\alpha s} l(X_s)\;ds=
E_\mu\int_0^\tau
 e^{-\alpha s}\Pi^\mu_s l\;ds
 $ can be proved in a similar way.
  \qed

Let again $(\bar\Omega,\bar\calf, (\bar\Pi_t)_{t\geq 0}, (\bar P_\nu)_{\nu\in \Delta_e},
(\bar\calf_t)_{t\geq 0})$ be the PDP on $\Delta_e$ defined in section \ref{PDPcanonical}.
We first need a technical result.

\begin{lemma}\label{trasportotempoarresto} For every $\mu\in\Delta$
and $\tau\in\calt^\mu$,  there
exists an $(\bar\calf^0_{t+})$-stopping time  $\bar\tau$ such that
$\tau(\omega)=\bar\tau(\Pi^\mu_.(\omega))$ for $P_\mu$-almost all $\omega\in\Omega$.
 $\bar\tau$ is also an $(\bar\calf_t)$-stopping time.
\end{lemma}

 \noindent {\bf Proof.}
 The last assertion of the lemma is obvious, since $\bar\calf^0_{t}\subset \bar\calf_t$
 and $(\bar\calf_{t})$ is a right-continuous filtration.

 Let us also note that the mapping $\Pi^\mu:\Omega\to\bar\Omega$,
 $\omega\mapsto \Pi^\mu_\cdot(\omega)$
 is measurable with respect to $\calf^0$ and $\bar\calf^0$: it is enough to note that
 for every $t$ the composition $\bar\Pi_t\circ \Pi^\mu: \Omega\to\Delta_e$ is
 measurable with respect to $\calf^0$ and $\calb(\Delta_e)$ since
 $(\bar\Pi_t\circ \Pi)( \omega)=\Pi_t^\mu(\omega)$.

 The rest of the proof consists in the construction of $\bar\tau$
 and will be divided into several steps.

 \noindent {\bf Step 1:} we show that $\caly^\mu_t:=\sigma(\caly^0_t,\caln^\mu)$
 coincides with $\calh:=\{A\subset\Omega\,:\, \exists B\in \caly^0_t,
 A \triangle B\in \caln^\mu \}$.

The argument is standard. The inclusion $\calh\subset \caly^\mu_t$ is easy and the inclusion
$\caly^\mu_t\subset \calh$ is proved by verifying that $\calh$ is a $\sigma$-algebra
containing $\caly^0_t$ and $\caln^\mu$.

 \noindent {\bf Step 2:}  for every $t\ge 0$ we have $\caly_t^0=(\Pi^\mu)^{-1}(\bar\calf_t^0)$.

 $\bar\calf_t^0$ is generated by the family of sets of the form
 \begin{equation}\label{genftzero}
    B=\{ \bar\omega\in\bar\Omega\;:\; \bar\omega(s_1)\in A_1,\ldots, \bar\omega(s_n)\in A_n\},
 \end{equation}
 for some integer $n$ and some $0\le s_1\le\ldots\le s_n\le t$, $A_i\in \calb(\Delta_e)$ ($1\le i\le n$).
 Then $(\Pi^\mu)^{-1}(B)$ is the set
 \begin{equation}\label{genftzerodue}
 \{ \omega\in\Omega\;:\; \Pi^\mu_{s_1}(\omega)\in A_1,\ldots, \Pi^\mu_{s_n}(\omega)\in A_n\},
 \end{equation}
 and this proves that $(\Pi^\mu)^{-1}(\bar\calf_t^0)\subset \caly_t^0$.

 Conversely, $\caly_t^0$ is generated by sets of the form (\ref{genftzerodue}), and each such
 set belongs to $(\Pi^\mu)^{-1}(\bar\calf_t^0)$ since it has the form
 $(\Pi^\mu)^{-1}(B)$ for $B$ given as in   (\ref{genftzero}). So we deduce that also
 $\caly_t^0\subset (\Pi^\mu)^{-1}(\bar\calf_t^0)$.

  \noindent {\bf Step 3:} we prove the result assuming that $\tau\in\calt^\mu$
  has a finite number of values.

  Excluding the trivial case where $\tau$ is constant, there exist an integer $n>1$,
  numbers $0\le t_1<\ldots<t_{n-1}<t_n=\infty$ and disjoint sets $A_i\subset \Omega$
  ($1\le i\le n$) such that $\cup_{i=1}^n A_i= \Omega$, $A_i\in \caly_{t_i}^\mu$ for
  $1\le i< n$ and $\tau(\omega)= \sum_{i=1}^n t_i1_{A_i}(\omega)$.
  By Step 1 there exist $N_i\in \caln^\mu $ and
$B_i\in \caly_{t_i}^0$ ($1\le i< n$) such that $A_i\triangle B_i=N_i$.
  By Step 2 there exist $C_i\in \bar\calf_{t_i}^0$ ($1\le i< n$) such that
  $(\Pi^\mu)^{-1}(C_i)=B_i$. Let us define $\bar\tau:\bar\Omega\to [0,\infty]$ setting
$$
\begin{array}{ll}
\bar\tau(\bar\omega)=t_1,
&
\bar\omega\in C_1,
\\
\bar\tau(\bar\omega)=t_i,
&
\bar\omega\in C_i\backslash (C_1\cup\ldots \cup C_{i-1}),\; 1<i<n,
\\
\bar\tau(\bar\omega)=\infty,
&
\bar\omega\notin C_1\cup\ldots \cup C_{n-1}.
\end{array}
$$
$\bar\tau$ is clearly a $(\bar\calf^0_t)$-stopping time.
To finish the proof of Step 3 it is enough to show that
\begin{equation}\label{composeccezz}
    \tau(\omega)=\bar\tau(\Pi^\mu(\omega)),
    \qquad
    \omega\notin N_1\cup\ldots \cup N_{n-1}.
\end{equation}
Let us fix $    \omega\notin N_1\cup\ldots \cup N_{n-1}$. First suppose that
$\Pi^\mu(\omega)\in C_i\backslash (C_1\cup\ldots \cup C_{i-1})$ for some $1<i<n$.
Then  in particular $\Pi^\mu(\omega)\in C_i$ and therefore $\omega\in B_i$ and
since $\omega\notin N_i$ it follows that $\omega\in A_i$ and
$\tau(\omega)=t_i=\bar\tau(\Pi^\mu(\omega))$. In the case
$\Pi^\mu(\omega)\in C_1$ the argument is similar. Finally if
$\Pi^\mu(\omega)\notin C_1\cup\ldots \cup C_{n-1}$ then for
 $1\le i<n$ we have $\omega\notin B_i$ and
since $\omega\notin N_i$ it follows that $\omega\notin A_i$, so we
deduce that $\tau(\omega)=\infty=\bar\tau(\Pi^\mu(\omega))$.
Now (\ref{composeccezz}) is proved.

  \noindent {\bf Step 4:} conclusion.

Given a general $\tau\in\calt^\mu$, let $\tau_n$ be a nonincreasing sequence of
$(\calf^\mu_{t})$-stopping times
such that $\tau_n\to\tau$
and each $\tau_n$ has a finite number of values. Let $\bar\tau_n$ be constructed
starting from $\tau_n$ as in Step 3, so that $P_\mu$-a.s. we have
$ \tau_n=\bar\tau_n(\Pi^\mu)$ for all $n$. Let us define
$\bar\tau(\bar\omega)=\liminf_{n\to\infty} \bar\tau_n(\bar\omega)$ for
all $\bar\omega\in \bar\Omega$ and let $A\subset \bar\Omega$ be the set where
$\lim_{n\to\infty} \bar\tau_n$ exists (finite or infinite). Then $\bar\tau$ is
a $(\bar\calf^0_{t+})$-stopping time. Moreover for  $P_\mu$-almost all $\omega\in\Omega$
we have $\bar\tau_n(\Pi^\mu_.(\omega))=\tau_n(\omega)\to \tau(\omega)$ which
shows that $\Pi^\mu_.(\omega)\in A$ and
$\bar\tau(\Pi^\mu_.(\omega))=\tau(\omega)$.
  \qed

\begin{lemma}\label{trasportocosto}
 For every $\mu\in\Delta$ and $\tau\in\calt^\mu$, let $\bar\tau$
denote an $(\bar\calf_t)$-stopping time such that
$\tau(\omega)=\bar\tau(\Pi^\mu_.(\omega))$ for $P_\mu$-almost all $\omega\in\Omega$.
Then we have
$$
 J(\mu,\tau)=\bar E_Q\left[e^{-\alpha \bar\tau}\bar\Pi_{\bar\tau} g+\int_0^{\bar\tau}
 e^{-\alpha s}\bar\Pi_s l\;ds\right],
 $$
 where $Q$ is the probability measure on $\Delta_e$ concentrated at points
 $H_a[\mu]$ ($a\in O$) such that
 $$
 Q(\{H_a[\mu]\})=\mu(h^{-1}(a)),\qquad a\in O.
 $$
\end{lemma}

 \noindent {\bf Proof.}
 Define a real function $\Phi$ on $\bar\Omega$ by
 $$
 \Phi(\bar\omega)=e^{-\alpha \bar\tau(\bar\omega)}\bar\Pi_{\bar\tau(\bar\omega)} (\bar\omega)g
 +\int_0^{\bar\tau(\bar\omega)}
 e^{-\alpha s}\bar\Pi_s (\bar\omega)l\;ds,
 \qquad \bar\omega\in\bar\Omega.
 $$
 It can be checked that $\Phi$ is bounded and $\bar\calf$-measurable.
 Taking $\bar\omega=\Pi^\mu_.(\omega)$ we obtain
 $\bar\Pi_s (\bar\omega)=\bar\omega(s)=\Pi^\mu_s(\omega)$ and
 for $\bar\tau(\bar\omega)<\infty$ we have
$ \bar\Pi_{\bar\tau(\bar\omega)} (\bar\omega)
 =\bar\omega(\bar\tau(\bar\omega))
 =\Pi^\mu_{\bar\tau(\bar\omega)}(\omega)$. Since $\bar\tau(\bar\omega)=
 \tau(\omega)$ $P_\mu$-a.s. we conclude that
 $$
 \Phi(\Pi^\mu_.(\omega))=
 e^{-\alpha \tau(\omega)}\Pi^\mu_{\tau(\omega)}(\omega) g+\int_0^{\tau(\omega)}
 e^{-\alpha s}\Pi^\mu_s(\omega) l\;ds,
 \qquad P_\mu-a.s.
 $$
 From lemma
\ref{passaggiotau}  and proposition \ref{leggefiltro}
we have
$$
 J(\mu,\tau)=\int_\Omega
  \Phi(\Pi^\mu_.(\omega)) \;P_\mu(d\omega)=
  \int_{\bar\Omega}
  \Phi(\bar\omega) \;\bar P_Q(d\bar\omega)=
  \bar E_Q\left[e^{-\alpha \bar\tau}\bar\Pi_{\bar\tau} g+\int_0^{\bar\tau}
 e^{-\alpha s}\bar\Pi_s l\;ds\right]. \qed
$$

We now formulate an auxiliary optimal stopping problem
with complete observation for the PDP
$(\bar\Omega,\bar\calf,(\bar P_\nu)_{\nu\in \Delta_e}, (\bar\Pi_t)_{t\geq 0},
(\bar\calf_t)_{t\geq 0})$. We denote by $\bar\calt$ the class of
$(\bar\calf_t)$-stopping times defined on $\bar\Omega$ and
we define the cost functional
$$
 \bar J(\nu,\bar\tau)=\bar E_\nu\left[e^{-\alpha \bar\tau}\bar\Pi_{\bar\tau} g+\int_0^{\bar\tau}
 e^{-\alpha s}\bar\Pi_s l\;ds\right], \qquad
\bar\tau\in\bar\calt,\;\nu\in\Delta_e,
 $$
 and the value function
 $$
 v(\nu) =\inf_{\bar\tau\in\bar\calt} \bar J(\nu,\bar\tau),
 \qquad \nu\in\Delta_e.
 $$
 Note that in this formulation the class of stopping times $\bar\calt$ does not
 depend on $\nu$, but it is sufficiently rich since it satisfies
 the usual conditions: compare remark \ref{filtrarresto}.

The proof of the following result can be found in \cite{Be}, Chapter VII, Theorem 6.1.
We do not repeat the assumptions of this theorem but we note that
they are all trivially verified in our situation; the only nontrivial verifications
are the properties of the semigroup $(R_t)$ which were proved
above in proposition \ref{pdpfeller}.

\begin{theorem}\label{arrottimopdp} The value function $v$ is continuous.
The random time defined by
$$
\bar\tau^*=\inf \{t\ge 0\,:\, \bar\Pi_t g=v(\bar\Pi_t)\}
$$
 belongs to $\bar\calt$ and it is optimal relatively to
every $\nu\in\Delta_e$:
$$
v(\nu)=\bar J(\nu,\bar\tau^{*}),\qquad
\nu\in\Delta_e.
$$
\end{theorem}

We are now ready to present the main result of this section.

\begin{theorem}\label{arrottimoossparz}
For $\mu\in\Delta$, let $(\Pi^\mu_t)$ be the filtering process
defined by equation (\ref{eqdefpifiltro}).
Define the random time
$$
\tau^{*,\mu}=\inf \{t\ge 0\,:\, \Pi^\mu_t g=v(\Pi^\mu_t)\}.
$$
Then $\tau^{*,\mu}$ belongs to $\calt^\mu$ and it is optimal relatively to
 $\mu$.
 Moreover the value function is given by the formula
$$
V(\mu)=\sum_{a\in O}\mu(h^{-1}(a))\,
v(H_a[\mu]),\qquad
\mu\in\Delta.
$$
\end{theorem}

\begin{remark}\begin{em}
\begin{enumerate}
\item[(i)]
The notation emphazises that $\tau^{*,\mu}$ depends on $\mu$, since
 $(\Pi^\mu_t)$  does.

\item[(ii)] If $\nu\in\Delta_e$ then
we have $V(\nu)=v(\nu)$. Indeed if $\nu\in\Delta_b$ for some $b\in O$ then
for all $a\in O$ we have
$\nu(h^{-1}(a))=1_{a=b}$ and $H_b[\nu]=\nu$.

\item[(iii)]
 $\tau^{*,\mu}$ is the first entry time of $\Pi^\mu$ in the contact set
 $\{\nu\in \Delta_e\,:\, \nu g= v(\nu)=V(\nu)\}$.
\end{enumerate}

\end{em}
  \end{remark}

 \noindent {\bf Proof.}
 Let $\mu\in\Delta$ be fixed. For arbitrary $\tau\in\calt^\mu$, by
 lemma \ref{trasportotempoarresto} there
exists an $(\bar\calf_t)$-stopping time  $\bar\tau$ such that
$\tau=\bar\tau(\Pi^\mu_.)$  $P_\mu$-a.s.
By lemma \ref{trasportocosto} and the definitions of $\bar P_Q$ and  $\bar J$,
\begin{equation}\label{trasffunzval}
    \begin{array}{lll}
 J(\mu,\tau)&=&\dis\bar E_Q\left[e^{-\alpha \bar\tau}\bar\Pi_{\bar\tau} g+\int_0^{\bar\tau}
 e^{-\alpha s}\bar\Pi_s l\;ds\right]
 \\
 &=&\dis
\sum_{a\in O} \mu(h^{-1}(a))
 \bar E_{H_a[\mu]}\left[e^{-\alpha \bar\tau}\bar\Pi_{\bar\tau} g+\int_0^{\bar\tau}
 e^{-\alpha s}\bar\Pi_s l\;ds\right]
 \\
 &=&\dis
\sum_{a\in O}\mu(h^{-1}(a))
 \bar J(H_a[\mu], \bar\tau).
 \end{array}
\end{equation}
 An application of
theorem \ref{arrottimopdp} shows that
\begin{equation}\label{trasffunzvaldue}
 J(\mu,\tau)
\ge \sum_{a\in O}\mu(h^{-1}(a))
 \bar J(H_a[\mu], \bar\tau^*)=
 \sum_{a\in O}\mu(h^{-1}(a))
v(H_a[\mu]).
\end{equation}

 Now consider $\tau^{*,\mu}$ and note that it is an element of $\calt^\mu$
 by the continuity of $v$. Moreover we have
$\tau^{*,\mu}=\bar\tau^*(\Pi^\mu_.)$  $P_\mu$-a.s. and so,
arguing as in (\ref{trasffunzval}) we obtain
$$ J(\mu,\tau^{*,\mu})=
\sum_{a\in O}\mu(h^{-1}(a))
 \bar J(H_a[\mu], \bar\tau^*).
$$
 Applying again
theorem \ref{arrottimopdp} we have
$$
 J(\mu,\tau^{*,\mu})=
 \sum_{a\in O}\mu(h^{-1}(a))
v(H_a[\mu]).
$$
Comparing with
(\ref{trasffunzvaldue})  we conclude that $\tau^{*,\mu}$ is optimal
relatively to $\mu$ and that the right-hand side of the last
formula equals $V(\mu)$.
  \qed

A satisfactory solution of the partially observed optimal stopping problem
should also include a characterization of its value function $V$.
The last assertion  of theorem \ref{arrottimoossparz} shows that
this problem is reduced to finding
 characterizations
of the value function $v$  of the optimal stopping problem for the PDP.
Since this is a fully observable problem
many results on analytical characterizations of  $v$ are known,
mostly in the form of obstacle problems.
For instance in
Theorem 6.1 of Chapter VII of
 \cite{Be}, mentioned above as theorem \ref{arrottimopdp}, it is also proved that
 $v$ satisfies the system of inequalities on $\Delta_e$
 \begin{equation}\label{inebellman}
    u\le \psi,
    \quad
    u\le e^{-\alpha t}R_{t} u+\int_0^{t}
 e^{-\alpha s}R_s L\;ds, \qquad t\ge 0,
 \end{equation}
 where $(R_t)$ is the transition semigroup
 of the PDP introduced before, the obstacle $\psi$ is simply $\psi(\nu)=\nu g$
 and the function $L$ is
 $L(\nu)=\nu l$ ($\nu\in \Delta_e$), and moreover $u$ is the maximum
 element among all real continuous functions on $\Delta_e$ satisfying
 (\ref{inebellman}).

In the specific case of optimal stopping for  PDPs many other
analytical characterizations of $v$ can be found in the literature,
(although sometimes under assumptions
slightly different from ours): see
 for instance \cite{Ga},
 \cite{Ga2} or  the monograph \cite{Da} and the references therein.


\end{document}